\documentclass[a4paper]{article}

\usepackage[cp1250]{inputenc}
\usepackage[IL2]{fontenc}

\usepackage{a4wide}
\usepackage[english]{babel}
\usepackage{euscript}
\usepackage{amstext,amsbsy,amscd,amssymb}
\usepackage{amsmath,enumitem}
\usepackage{amsfonts}
\usepackage{graphics}
\usepackage{graphicx}
\usepackage{microtype}
\usepackage{color}
\include{diagxy}
\allowdisplaybreaks

\let\rarr=\rightarrow

\let\larr=\leftarrow

\let\veps=\varepsilon
\let\mcal=\mathcal
\let\mfrak=\mathfrak
\let\eus=\EuScript

\def\N{\mathbb{N}}
\def\Z{\mathbb{Z}}
\def\R{\mathbb{R}}
\def\C{\mathbb{C}}

\def\Mod{\mathop {\rm Mod} \nolimits}

\def\End{\mathop {\rm End} \nolimits}
\def\Ker{\mathop {\rm Ker} \nolimits}

\def\Hom{\mathop {\rm Hom} \nolimits}

\def\ad{\mathop {\rm ad} \nolimits}
\def\gr{\mathop {\rm gr} \nolimits}

\def\GL{\mathop {\rm GL} \nolimits}

\def\id{{\rm id}}

\def\Ad{\mathop {\rm Ad} \nolimits}

\def\Ind{\mathop {\rm Ind} \nolimits}

\def\htt{\mathop {\rm ht} \nolimits}

\def\Cas{\mathop {\rm Cas} \nolimits}

\def\Specm{\mathop {\rm Specm} \nolimits}

\def\GInd{\mathop {\rm GInd} \nolimits}

\newsymbol\squares 1003

\long\def\proof #1{\noindent \emph{Proof.}\ #1 \hfill $\squares$

\medskip}

\newcounter{num}[section]
\numberwithin{equation}{section}
\numberwithin{num}{section}

\long\def\definition #1 {\refstepcounter{num} \noindent {\bf
Definition \thenum.} #1

\medskip}

\long\def\theorem #1{\refstepcounter{num} \noindent {\bf Theorem
\thenum.} #1

\medskip}

\long\def\lemma #1{\refstepcounter{num}  \noindent {\bf Lemma
\thenum.} #1

\medskip}

\long\def\proposition #1{\refstepcounter{num}  \noindent {\bf Proposition
\thenum.} #1

\medskip}

\long\def\remark #1{\noindent {\bf Remark.}\ #1

\medskip}

\newcommand*\riso{%
  \xrightarrow[]{\raisebox{-0.25em}{\smash{\ensuremath{\sim}}}}%
}

\makeatletter
\newcommand*\if@single[3]{%
  \setbox0\hbox{${\mathaccent"0362{#1}}^H$}%
  \setbox2\hbox{${\mathaccent"0362{\kern0pt#1}}^H$}%
  \ifdim\ht0=\ht2 #3\else #2\fi
  }
\newcommand*\rel@kern[1]{\kern#1\dimexpr\macc@kerna}
\newcommand*\widebar[1]{\@ifnextchar^{{\wide@bar{#1}{0}}}{\wide@bar{#1}{1}}}
\newcommand*\wide@bar[2]{\if@single{#1}{\wide@bar@{#1}{#2}{1}}{\wide@bar@{#1}{#2}{2}}}
\newcommand*\wide@bar@[3]{%
  \begingroup
  \def\mathaccent##1##2{%
    \if#32 \let\macc@nucleus\first@char \fi
    \setbox\z@\hbox{$\macc@style{\macc@nucleus}_{}$}%
    \setbox\tw@\hbox{$\macc@style{\macc@nucleus}{}_{}$}%
    \dimen@\wd\tw@
    \advance\dimen@-\wd\z@
    \divide\dimen@ 3
    \@tempdima\wd\tw@
    \advance\@tempdima-\scriptspace
    \divide\@tempdima 10
    \advance\dimen@-\@tempdima
    \ifdim\dimen@>\z@ \dimen@0pt\fi
    \rel@kern{0.6}\kern-\dimen@
    \if#31
      \overline{\rel@kern{-0.6}\kern\dimen@\macc@nucleus\rel@kern{0.4}\kern\dimen@}%
      \advance\dimen@0.4\dimexpr\macc@kerna
      \let\final@kern#2%
      \ifdim\dimen@<\z@ \let\final@kern1\fi
      \if\final@kern1 \kern-\dimen@\fi
    \else
      \overline{\rel@kern{-0.6}\kern\dimen@#1}%
    \fi
  }%
  \macc@depth\@ne
  \let\math@bgroup\@empty \let\math@egroup\macc@set@skewchar
  \mathsurround\z@ \frozen@everymath{\mathgroup\macc@group\relax}%
  \macc@set@skewchar\relax
  \let\mathaccentV\macc@nested@a
  \if#31
    \macc@nested@a\relax111{#1}%
  \else
    \def\gobble@till@marker##1\endmarker{}%
    \futurelet\first@char\gobble@till@marker#1\endmarker
    \ifcat\noexpand\first@char A\else
      \def\first@char{}%
    \fi
    \macc@nested@a\relax111{\first@char}%
  \fi
  \endgroup
}
\makeatother

\newcommand\rsmraise[1]{%
  \ifx#1\displaystyle .8\else
    \ifx#1\textstyle .8\else
      \ifx#1\scriptstyle .6\else
        .45%
      \fi
    \fi
  \fi}

\title{Twisting functors and Gelfand--Tsetlin modules \\
over semisimple Lie algebras}

\author{Vyacheslav Futorny, Libor Křižka}

\AtEndDocument{\bigskip{\footnotesize%
  (V.\,Futorny) \textsc{Instituto de Matemática e Estatística, Universidade de Sa\~{o} Paulo, Caixa Postal 66281,\\ Sa\~{o} Paulo, CEP 05315-970, Brasil} \par
  \textit{E-mail address}: \texttt{futorny@ime.usp.br} \par
  \addvspace{\medskipamount}
  (L.\,Křižka) \textsc{Instituto de Matemática e Estatística, Universidade de Sa\~{o} Paulo, Caixa Postal 66281,\\ Sa\~{o} Paulo, CEP 05315-970, Brasil} \par
  \textit{E-mail address}: \texttt{krizka.libor@gmail.com} \par
}}

\date{\today}

\begin{document}

\maketitle

\begin{abstract}
We associate to an arbitrary positive root $\alpha$ of a complex semisimple finite-dimensional Lie algebra $\mfrak{g}$ a twisting endofunctor $T_\alpha$ of the category  of $\mfrak{g}$-modules. We apply this functor to generalized Verma modules in the category $\mcal{O}(\mfrak{g})$ and construct a family of $\alpha$-Gelfand--Tsetlin modules with finite $\Gamma_\alpha$-multiplicities, where $\Gamma_{\alpha}$ is a commutative $\C$-subalgebra of the universal enveloping algebra of $\mfrak{g}$ generated by a Cartan subalgebra of $\mfrak{g}$ and by the Casimir element of the $\mfrak{sl}(2)$-subalgebra corresponding to the root $\alpha$. This covers classical results of Andersen and Stroppel when $\alpha$ is a simple root and previous results of the authors in the case when $\mfrak{g}$ is a complex simple Lie algebra and  $\alpha$ is the maximal root of $\mfrak{g}$. The significance of constructed modules is that they are Gelfand--Tsetlin modules with respect to any commutative $\C$-subalgebra of the universal enveloping algebra of $\mfrak{g}$ containing $\Gamma_\alpha$. Using the Beilinson--Bernstein correspondence we give a geometric realization of these modules together with their explicit description. We also identify a tensor subcategory of the category of $\alpha$-Gelfand--Tsetlin modules which contains constructed modules as well as the category $\mcal{O}(\mfrak{g})$.

\medskip

\medskip

\medskip

\end{abstract}

\thispagestyle{empty}

\tableofcontents


\section*{Introduction}
\addcontentsline{toc}{section}{Introduction}

The theory of \emph{Gelfand--Tsetlin modules} has been developed extensively in recent years following a number of important results obtained  in \cite{Futorny-Grantcharov-Ramirez2016}, \cite{Futorny-Grantcharov-Ramirez2017}, \cite{Futorny-Ramirez-Zhang2016}, \cite{Futorny-Grantcharov-Ramirez2017a}, \cite{Zadunaisky2017}, \cite{Vishnyakova2018}, \cite{Vishnyakova2017}, \cite{Ramirez-Zadunaisky2018}, \cite{Futorny-Krizka2017}, \cite{Early-Mazorchuk-Vishnyakova2017}, \cite{Hartwig2017}, \cite{Futorny-Grantcharov-Ramirez-Zadunaisky2018}, \cite{Mazorchuk-Vishnyakova2018}, \cite{Futorny-Grantcharov-Ramirez-Zadunaisky2018a}. These modules are connected with the study of Gelfand--Tsetlin integrable systems (e.g.\ \cite{Guillemin-Sternberg1983}, \cite{Kostant-Wallach2006}, \cite{Graev2007}), BGG differential operators (e.g.\ \cite{Early-Mazorchuk-Vishnyakova2017}, \cite{Futorny-Grantcharov-Ramirez-Zadunaisky2018}), tensor product categorifications and KLRW-algebras in \cite{Kamnitzer-Tingley-Webster-Weekes-Yacobi2018}.

Gelfand--Tsetlin modules were originally defined and studied in \cite{Drozd-Ovsienko-Futorny1991} and \cite{Drozd-Futorny-Ovsienko1994} aiming to parameterize simple modules for the universal enveloping algebra of a complex simple Lie algebra by characters of a certain maximal commutative $\C$-subalgebra -- the \emph{Gelfand--Tsetlin subalgebra}. Up to now the theory is mainly evolved for simple Lie algebras of type $A$, that is $U(\mfrak{sl}(n))$ and its generalizations. In \cite{Futorny-Krizka2017} a new approach was developed which allows us to construct certain Gelfand--Tsetlin-like modules for an arbitrary simple finite-dimensional Lie algebra $\mfrak{g}$.

Let $\theta$ be the maximal root of $\mfrak{g}$ and let $\mfrak{s}_\theta$ be the Lie subalgebra of $\mfrak{g}$ based on the root $\theta$, and hence isomorphic to $\mfrak{sl}(2)$. Let us denote by $\Gamma_\theta$ the commutative $\C$-subalgebra of the universal enveloping algebra $U(\mfrak{g})$ of $\mfrak{g}$ generated by a Cartan subalgebra of $\mfrak{g}$ and by the center of $U(\mfrak{s}_\theta)$. The \emph{$\theta$-Gelfand--Tsetlin modules}, modules admitting a locally finite action of $\Gamma_\theta$, were studied in \cite{Futorny-Krizka2017}.
In the case of $\mfrak{sl}(n)$ such modules give examples of the so-called \emph{partial Gelfand--Tsetlin modules} studied in \cite{Futorny-Ovsienko-Saorin2011}. In particular, a family of  $\theta$-Gelfand--Tsetlin modules $W^\mfrak{g}_\mfrak{b}(\lambda, \theta)$ with finite $\Gamma_\theta$-multiplicities were constructed explicitly using a free field realization.
Here $\mfrak{b}$ is a Borel subalgebra of $\mfrak{g}$ with the corresponding Cartan subalgebra $\mfrak{h}$ and $\lambda\in \mfrak{h}^*$. Using the Beilinson--Bernstein correspondence one can view $W^\mfrak{g}_\mfrak{b}(\lambda, \theta)$ as the vector space of `$\delta$-functions' on the flag variety $G/B$ supported at the $1$-dimensional subvariety being an orbit of a unipotent subgroup of $G$ and going through the point $eB$.

Let $\Gamma$ be a commutative $\C$-subalgebra of the universal enveloping algebra $U(\mfrak{g})$ of $\mfrak{g}$. Then we denote by $\mcal{H}(\mfrak{g},\Gamma)$ the category of $\Gamma$-weight $\mfrak{g}$-modules and by $\mcal{H}_{\rm fin}(\mfrak{g},\Gamma)$ the full subcategory of $\mcal{H}(\mfrak{g},\Gamma)$ consisting of $\Gamma$-weight $\mfrak{g}$-modules with finite-dimensional $\Gamma$-weight subspaces. The modules $W^\mfrak{g}_\mfrak{b}(\lambda, \theta)$ belong to the category $ \mcal{H}_{\rm fin}(\mfrak{g},\Gamma_{\theta})$ and as an easy consequence  also to the category $\mcal{H}_{\rm fin}(\mfrak{g},\Gamma)$ for any commutative $\C$-subalgebra $\Gamma$ of $U(\mfrak{g})$ containing $\Gamma_\theta$.

The purpose of the current paper is to extend the construction of $W^\mfrak{g}_\mfrak{b}(\lambda, \theta)$ to any positive root $\alpha$ for an arbitrary complex semisimple finite-dimensional Lie algebra. We show that such modules can be obtained from the corresponding Verma modules by applying the twisting functor. This functor was introduced by Deodhar \cite{Deodhar1980} and
was used successfully by Mathieu \cite{Mathieu2000} to classify simple torsion free weight modules with finite-dimensional weight spaces over complex simple Lie algebras of type $A$ and $C$. If $\alpha$ is a simple root, then the twisting functor is well understood, it is related to the Arkhipov's twisting functor \cite{Arkhipov2004} on the category $\mcal{O}$. We are mostly interested in the case when $\alpha$ is not simple root. In addition, we replace the Borel subalgebra $\mfrak{b}$ by any standard parabolic subalgebra $\mfrak{p}$ of $\mfrak{g}$. The corresponding \emph{$\alpha$-Gelfand--Tsetlin} module $W^\mfrak{g}_\mfrak{p}(\lambda, \alpha)$ is obtained by applying the twisting functor to the generalized Verma module $M^\mfrak{g}_\mfrak{p}(\lambda)$ induced from the finite-dimensional simple $\mfrak{p}$-modules with highest weight $\lambda$. We show that all modules $W^\mfrak{g}_\mfrak{p}(\lambda, \alpha)$ are cyclic weight modules with respect to the Cartan subalgebra $\mfrak{h}$. Moreover, every element of $\Gamma_\alpha$ has a Jordan decomposition on $W^\mfrak{g}_\mfrak{p}(\lambda, \alpha)$ with Jordan cells of size at most $2$. The same holds for the elements of any commutative $\C$-subalgebra $\Gamma$ of $U(\mfrak{g})$ containing $\Gamma_\alpha$. In the case $\mfrak{g}=\mfrak{sl}(3)$ the subalgebra $\Gamma_\theta$ is generically diagonalizable on $W^\mfrak{g}_\mfrak{b}(\lambda, \theta)$ and also the latter module is generically simple, $\theta$ is the maximal root of $\mfrak{g}$, see \cite{Futorny-Krizka2017}. We expect the same properties of modules $W^\mfrak{g}_\mfrak{p}(\lambda, \alpha)$ in general. This will be addressed in a subsequent paper.

Besides, we provide the $\mcal{D}$-module realization of modules $W^\mfrak{g}_\mfrak{b}(\lambda, \alpha)$ using the Beilinson--Bernstein correspondence (Theorem \ref{theo-realiz}). We also give explicit formulas for the Lie algebra action in this realization (Theorem \ref{thm:Weyl realization}).

\begin{table}[ht]
\centering
\renewcommand{\arraystretch}{1.5}
\begin{tabular}{c|c|c}
  generalized Verma module & Gelfand--Tsetlin module & Gelfand--Tsetlin module\\
  $M^\mfrak{g}_\mfrak{p}(\lambda)$, $\alpha \in \Delta^+_\mfrak{u}$ & $W^\mfrak{g}_\mfrak{p}(\lambda,\alpha)$, $\alpha \in \Delta^+_\mfrak{u} \cap \Pi$ & $W^\mfrak{g}_\mfrak{p}(\lambda,\alpha)$, $\alpha \in \Delta^+_\mfrak{u} \setminus \Pi$ \\
    \hline
  central character $\chi_{\lambda+\rho}$ & central character $\chi_{\lambda+\rho}$ & central character $\chi_{\lambda+\rho}$ \\
  cyclic module & cyclic module & cyclic module \\[-2mm]
  \parbox[c]{4cm}{\begin{center} weight $\mfrak{g}$-module with \\ finite-dimensional \\  weight spaces \end{center}} &  \parbox[c]{4cm}{\begin{center} weight $\mfrak{g}$-module with \\ finite-dimensional \\  weight spaces \end{center}} &  \parbox[c]{4cm}{\begin{center} weight $\mfrak{g}$-module with \\ infinite-dimensional \\  weight spaces \end{center}} \\[-4mm]
  \parbox[c]{4cm}{\begin{center}  $\Gamma_\alpha$-weight $\mfrak{g}$-module with \\ finite $\Gamma_\alpha$-multiplicities \\   \end{center}} &  \parbox[c]{4cm}{\begin{center}  $\Gamma_\alpha$-weight $\mfrak{g}$-module  with \\ finite $\Gamma_\alpha$-multiplicities  \end{center}} &  \parbox[c]{4cm}{\begin{center}  $\Gamma_\alpha$-weight $\mfrak{g}$-module  with \\ finite $\Gamma_\alpha$-multiplicities  \end{center}}
\end{tabular}
\caption{Comparison of generalized Verma modules and Gelfand--Tsetlin modules}  \label{tab:comparision}
\end{table}

One of the  inspirations and motivations for our original quest was the question of A.\,Kleshchev \cite{Kleshchev2018} about the existence of tensor categories of $\Gamma$-Gelfand--Tsetlin modules beyond the category $\mcal{O}$ for $\mfrak{g}=\mfrak{gl}(n)$ and $\mfrak{g}=\mfrak{sl}(n)$, where $\Gamma$ is the Gelfand--Tsetlin subalgebra of $\mfrak{g}$. We observe that $\Gamma_\alpha$-Gelfand--Tsetlin modules form a tensor category which contains the category $\mcal{O}$. In particular, if $\mfrak{g}=\mfrak{gl}(n)$ then $\Gamma$-Gelfand--Tsetlin modules which are $\Gamma_\alpha$-Gelfand--Tsetlin modules form a tensor category if $\Gamma$ contains $\Gamma_\alpha$. We provide a recipe how to construct more general tensor categories of partial Gelfand--Tsetlin modules.

Let us briefly summarize the content of our paper. In Section \ref{sec:Twisting functors and Gelfand--Tsetlin} we recall a general notion of twisting functors for rings and introduce $\Gamma$-Gelfand--Tsetlin modules for a commutative $\C$-subalgebra $\Gamma$ of the universal enveloping algebra $U(\mfrak{g})$ of a Lie algebra $\mfrak{g}$. Section \ref{sec:Gelfand--Tsetlin} is devoted to studying of $\alpha$-Gelfand--Tsetlin modules over complex semisimple Lie algebras. We also introduce various tensor categories of $\alpha$-Gelfand--Tsetlin modules. In Section \ref{sec:Twisting functor and Verma modules}, we define the twisting functor $T_\alpha$ for any positive root $\alpha$ of a complex semisimple Lie algebra $\mfrak{g}$ as an endofunctor of the category $\mcal{M}(\mfrak{g})$ of $\mfrak{g}$-modules and describe main properties of this functor. By applying the twisting functor $T_\alpha$ on generalized Verma modules $M^\mfrak{g}_\mfrak{p}(\lambda)$ we construct $\alpha$-Gelfand--Tsetlin modules $W^\mfrak{g}_\mfrak{p}(\lambda,\alpha)$ with finite $\Gamma_\alpha$-multiplicities (Theorem \ref{thm:GT modules}). Basic characteristics of these modules together with their geometric  realization are described in Section \ref{sec:Realizations}.


\section{Preliminaries}
\label{sec:Twisting functors and Gelfand--Tsetlin}

We denote by $\C$, $\R$, $\Z$, $\N_0$ and $\N$ the set of complex numbers, real numbers, integers, non-negative integers and positive integers, respectively. All algebras and modules are considered over the field of complex numbers.


\subsection{Twisting functors}

We recall the definition and basis properties of twisting functors for rings.
\medskip

\definition{Let $R$ be a ring. A multiplicative set $S$ in $R$ is called a left denominator set if it satisfies the following two conditions
\begin{enumerate}[topsep=0pt,itemsep=0pt,parsep=0pt]
  \item[i)] $Sr \cap Rs \neq \emptyset$ for all $r \in R$ and $s \in S$ (the set $S$ is left permutable);
  \item[ii)] if $rs=0$ for $r \in R$ and $s \in S$, then there exists $s' \in S$ such that $s'r=0$ (the set $S$ is left reversible).
\end{enumerate}
The above two conditions are usually called Ore's conditions.}

If $S$ is a left denominator set in a ring $R$, then we may construct a left ring of fractions for $R$ with respect to $S$. Let us recall that a left ring of fractions for $R$ with respect to $S$ is a ring homomorphism $f \colon R \rarr T$ such that
\begin{enumerate}[topsep=0pt,itemsep=0pt,parsep=0pt]
 \item[i)] $f(s)$ is a unit in $T$ for all $s \in S$;
 \item[ii)] any element of $T$ can be written in the form $f(s)^{-1}f(r)$ for some $s \in S$ and $r \in R$;
 \item[iii)] $\ker f = \{r \in R;\, sr=0\ \text{for some $s \in S$}\}$.
\end{enumerate}
Let us note that a left ring of fractions is unique up to isomorphism, we will denote it by $S^{-1}R$.

Further, if $M$ is a left $R$-module, then a module of fractions for $M$ with respect to $S$ is an $R$-module homomorphism $\varphi \colon M \rarr N$, where $N$ is a left $S^{-1}R$-module, such that
\begin{enumerate}[topsep=0pt,itemsep=0pt,parsep=0pt]
  \item[i)] any element of $N$ can be written in the form $s^{-1}\varphi(a)$ for some $s \in S$ and $a \in M$;
  \item[ii)] $\ker \varphi = \{a \in M;\, sa=0\ \text{for some $s \in S$}\}$.
\end{enumerate}
Let us note that a module of fractions is unique up to isomorphism, we will denote it by $S^{-1}M$.
\medskip

\theorem{Let $S$ be a left denominator set in a ring $R$ and let $M$ be a left $R$-module. Then we have
\begin{align}
  S^{-1}R \otimes_R M \simeq S^{-1}M,
\end{align}
where the mapping is given by $s^{-1}r \otimes a \mapsto s^{-1}ra$ for $s \in S$, $r \in R$ and $a \in M$.}

\lemma{Let $R$ be a ring and let $t \in R$ be a locally $\ad$-nilpotent regular element. Then the multiplicative set  $S_t=\{t^n;\, n\in \N_0\}$ in $R$ is a left denominator set.}

\proof{To check Ore's conditions for $S_t$, it is sufficient to verify them for the generator $t$ of $S_t$. As $t$ is a locally $\ad$-nilpotent element, for any $r \in R$ we have $\ad(t)^n(r) = 0$ for some $n \in \N$. Then the identity $\ad(t)^n(r) = \sum_{k=0}^n (-1)^k\binom{n}{k} t^{n-k}rt^k=0$ can be written as $t^nr=r't$ for some $r' \in R$. Moreover, since $t$ is a regular element, Ore's conditions are satisfied for $t$.}

Let $R$ be a ring and let $t \in R$ be a locally $\ad$-nilpotent regular element. A left ring of fractions for $R$ with respect to $S_t$ we will denote by $R_{(t)}$, and similarly a module of fractions for $M$ with respect to $S_t$ we will denote by $M_{(t)}$. Since $t$ is a regular element of $R$, the ring homomorphism $R \rarr R_{(t)}$ is injective. Hence, we may regard $R$ as a subring of $R_{(t)}$. We also have $(R,R)$-bimodule $R_{(t)}/R$, which enables us to define the twisting functor
\begin{align*}
  T_t \colon \Mod(R) \rarr \Mod(R)
\end{align*}
by
\begin{align}
  T_t(M) = (R_{(t)}/R) \otimes_R M
\end{align}
for $M \in \Mod(R)$. Besides, we introduce the $1$-parameter family of automorphisms $\theta_t^\nu \colon R_{(t)} \rarr R_{(t)}$ by
\begin{align}
  \theta_t^\nu(r) = t^\nu r t^{-\nu} = \sum_{k=0}^\infty \binom{\nu+k-1}{k} t^{-k}\ad(t)^k(r)
\end{align}
for $\nu \in \C$ and $r \in R_{(t)}$. Let us note that the sum on the right hand side is well defined since $t$ is a locally $\ad$-nilpotent element of $R$. The following lemma is straightforward.
\medskip

\lemma{Let $R$ be a ring and $t \in R$ a locally $\ad$-nilpotent regular element. Then we have
$\theta_t^\mu \circ \theta_t^\nu = \theta_t^{\mu+\nu}$ for $\mu, \nu \in \C$.}

\proposition{Let $R$ be a ring and $t \in R$ a locally $\ad$-nilpotent regular element. Then the twisting functor $T_t$ is right exact.}

In this article, we will mostly interested in two special cases when the ring $R$ is the universal enveloping algebra of a complex semisimple Lie algebra or the Weyl algebra of a complex vector space.


\subsection{$\Gamma$-Gelfand--Tsetlin modules}

Let $\mfrak{g}$ be a complex semisimple finite-dimensional Lie algebra and let $\mfrak{h}$ be a Cartan subalgebra of $\mfrak{g}$. We denote by $\Delta$ the root system of $\mfrak{g}$ with respect to $\mfrak{h}$, by $\Delta^+$ a positive root system in $\Delta$ and by $\Pi \subset \Delta^+$ the set of simple roots. For $\alpha \in \Delta^+$, let $h_\alpha \in \mfrak{h}$ be the corresponding coroot and let $e_\alpha$ and $f_\alpha$ be basis of $\mfrak{g}_\alpha$ and $\mfrak{g}_{-\alpha}$, respectively, defined by the requirement $[e_\alpha, f_\alpha] = h_\alpha$. We also set
\begin{align*}
  Q = \sum_{\alpha \in \Delta^+} \Z \alpha \qquad \text{and} \qquad Q_+ = \sum_{\alpha \in \Delta^+} \N_0 \alpha
\end{align*}
together with
\begin{align*}
  \Lambda = \{\lambda \in \mfrak{h}^*;\, (\forall \alpha \in \Pi)\, \lambda(h_\alpha) \in \Z\} \qquad \text{and} \qquad \Lambda^+ = \{\lambda \in \mfrak{h}^*;\, (\forall \alpha \in \Pi)\, \lambda(h_\alpha) \in \N_0\}.
\end{align*}
We call $Q$ the root lattice and $\Lambda$ the weight lattice. Further, we define the Weyl vector $\rho \in \mfrak{h}^*$ by
\begin{align*}
 \rho =  {1 \over 2} \sum_{\alpha \in \Delta^+} \alpha.
\end{align*}
The standard Borel subalgebra $\mfrak{b}$ of $\mfrak{g}$ is defined through $\mfrak{b} = \mfrak{h} \oplus \mfrak{n}$ with the nilradical $\mfrak{n}$ and the opposite nilradical $\widebar{\mfrak{n}}$ given by
\begin{align*}
  \mfrak{n} = \bigoplus_{\alpha \in \Delta^+} \mfrak{g}_\alpha \qquad \text{and} \qquad \widebar{\mfrak{n}} = \bigoplus_{\alpha \in \Delta^+} \mfrak{g}_{-\alpha}.
\end{align*}
Besides, we have the corresponding triangular decomposition
\begin{align*}
  \mfrak{g} = \widebar{\mfrak{n}} \oplus \mfrak{h} \oplus \mfrak{n}
\end{align*}
of the Lie algebra $\mfrak{g}$. Further, let $(\cdot\,,\cdot)_\mfrak{g} \colon \mfrak{g} \otimes_\C \mfrak{g} \rarr \C$ be the Cartan--Killing form on $\mfrak{g}$. Whenever $\alpha \in \mfrak{h}^*$ satisfies  $(\alpha,\alpha)_\mfrak{g} \neq 0$, we define $s_\alpha \in \GL(\mfrak{h}^*)$ by
\begin{align*}
  s_\alpha(\gamma) = \gamma - {2(\alpha,\gamma)_\mfrak{g} \over (\alpha,\alpha)_\mfrak{g}}\, \alpha
\end{align*}
for $\gamma \in \mfrak{h}^*$. The subgroup $W$ of $\GL(\mfrak{h}^*)$ defined through
\begin{align*}
  W = \langle s_\alpha;\, \alpha \in \Pi \rangle
\end{align*}
is called the Weyl group of $\mfrak{g}$. Let us note that $W$ is a finite Coxeter group. Let $\mfrak{z}_\mfrak{g}$ be the center of $U(\mfrak{g})$ and let $\gamma \colon \mfrak{z}_\mfrak{g} \rarr U(\mfrak{h})$ be the Harish-Chandra homomorphism (its image coincides with the set of $W$-invariant elements of $U(\mfrak{h})$). For each $\lambda \in \mfrak{h}^*$ we define a central character $\chi_\lambda \colon \mfrak{z}_\mfrak{g} \rarr \C$ by
\begin{align*}
  \chi_\lambda(z)=(\gamma(z))(\lambda)
\end{align*}
for $z \in \mfrak{z}_\mfrak{g}$, where we identify $U(\mfrak{h})$ with the $\C$-algebra of polynomial functions on $\mfrak{h}^*$.
\medskip

For a commutative $\C$-algebra $\Gamma$ we denote by $\Hom(\Gamma,\C)$ the set of all characters of $\Gamma$, i.e.\ $\C$-algebra homomorphisms from $\Gamma$ to $\C$. Let us note that if $\Gamma$ is finitely generated, then there is a natural identification between the set $\Hom(\Gamma,\C)$ of all characters of $\Gamma$ and the set $\Specm \Gamma$ of all maximal ideals of $\Gamma$, which corresponds to a complex algebraic variety. Let $M$ be a $\Gamma$-module. For each $\chi \in \Hom(\Gamma,\C)$  we set
\begin{align}
  M_\chi = \{v \in M;\, (\exists k \in \N)\,(\forall a \in \Gamma)\, (a-\chi(a))^kv=0\}.
\end{align}
When $M_\chi \neq\{0\}$, we say that $\chi$ is a $\Gamma$-weight of $M$, the vector space $M_\chi$ is called the $\Gamma$-weight subspace of $M$ with weight $\chi$ and the elements of $M_\chi$ are $\Gamma$-weight vectors with weight $\chi$. If a $\Gamma$-module $M$ satisfies
\begin{align}
  M = \bigoplus_{\chi \in \Hom(\Gamma,\C)} M_\chi,
\end{align}
then we call $M$ a \emph{$\Gamma$-weight module}. The dimension of the vector space $M_\chi$ will be called the \emph{$\Gamma$-multiplicity} of $\chi$ in $M$. Let us note that any submodule or factor-module of a $\Gamma$-weight module is also a $\Gamma$-weight module.

Let $\Gamma$ be a commutative $\C$-subalgebra of the universal enveloping algebra $U(\mfrak{g})$ of $\mfrak{g}$. Then we denote by $\mcal{H}(\mfrak{g},\Gamma)$ the category of all $\Gamma$-weight $\mfrak{g}$-modules and by $\mcal{H}_{\rm fin}(\mfrak{g},\Gamma)$ the full subcategory of $\mcal{H}(\mfrak{g},\Gamma)$ consisting of $\Gamma$-weight $\mfrak{g}$-modules with finite-dimensional $\Gamma$-weight subspaces. If $\Gamma$ and $\Gamma'$ are commutative $\C$-subalgebras of $U(\mfrak{g})$ satisfying $\Gamma \subset \Gamma'$, then we have the following obvious inclusions
\begin{align}
 \mcal{H}_{\rm fin}(\mfrak{g},\Gamma) \subset \mcal{H}_{\rm fin}(\mfrak{g},\Gamma') \subset \mcal{H}(\mfrak{g},\Gamma') \subset \mcal{H}(\mfrak{g},\Gamma),
\end{align}
which are strict in general.
If $\Gamma$ contains the Cartan subalgebra $\mfrak{h}$, then a $\Gamma$-weight $\mfrak{g}$-module $M$ is called a \emph{$\Gamma$-Gelfand--Tsetlin} module. In particular, a usual \emph{weight} $\mfrak{g}$-module $M$, i.e.\ a $\mfrak{g}$-module $M$ satisfying
\begin{align*}
M=\bigoplus_{\lambda\in \mfrak{h}^*}\! M_\lambda,
\end{align*}
where $M_\lambda = \{v \in M;\, (\forall h \in \mfrak{h})\, hv=\lambda(h)v\}$, is a $U(\mfrak{h})$-Gelfand--Tsetlin module.
\medskip

Let $\{e,h,f\}$ denotes the standard basis of the Lie algebra $\mfrak{sl}(2)$. The vector subspace $\mfrak{h}=\C h$ is a Cartan subalgebra of $\mfrak{sl}(2)$. If we define $\alpha \in \mfrak{h}^*$ by $\alpha(h)=2$, then the root system of $\mfrak{sl}(2)$ with respect to $\mfrak{h}$ is $\Delta=\{\pm\alpha\}$ with $\Delta^+=\{\alpha\}$. The standard Borel subalgebra $\mfrak{b}$ of $\mfrak{sl}(2)$ is defined as $\mfrak{b}=\C h \oplus \C e$ with the nilradical $\mfrak{n}=\C e$ and the opposite nilradical $\widebar{\mfrak{n}} = \C f$. Besides, the center $\mfrak{z}_{\mfrak{sl}(2)}$ of $U(\mfrak{sl}(2))$ is freely generated by the quadratic Casimir element $\Cas$ given by
\begin{align*}
  \Cas = ef+fe+{\textstyle {1 \over 2}} h^2.
\end{align*}
If $M$ is an $\mfrak{sl}(2)$-module then $M$ is also a $\mfrak{z}_{\mfrak{sl}(2)}$-module. The following proposition is standard but we include the details for convenience of the reader.
\medskip

\proposition{\label{prop:Jordan blocks}
\begin{enumerate}[topsep=3pt,itemsep=0pt]
\item[i)] Let $M$ be an $\mfrak{sl}(2)$-module. If for a weight vector $v \in M$ with weight $\lambda \in \mfrak{h}^*$ there exists $n \in \N$ such that $f^n v = 0$, then we have
    \begin{align}
      \prod_{k=0}^{n-1} (z-\chi_{\lambda-\rho-k\alpha}(z))v = 0 \label{eq:Casimir equation}
    \end{align}
    for all $z \in \mfrak{z}_{\mfrak{sl}(2)}$.
\item[ii)] Let $M$ be a weight $\mfrak{sl}(2)$-module which is locally $\widebar{\mfrak{n}}$-finite. Then $M$ is a $\mfrak{z}_{\mfrak{sl}(2)}$-weight module. Moreover, the $\mfrak{z}_{\mfrak{sl}(2)}$-weights of $M$ may be only of the form $\chi_{\lambda-\rho}$ for $\lambda \in \mfrak{h}^*$ and Jordan blocks corresponding to the $\mfrak{z}_{\mfrak{sl}(2)}$-weight $\chi_{\lambda-\rho}$ are of size at most $2$ if $w(\lambda-\rho)-\rho \in \Lambda^+$ for some $w \in W$ and of size $1$ otherwise.
\end{enumerate}
}

\proof{i) It is sufficient to prove the statement for the generator $\Cas$ of $\mfrak{z}_{\mfrak{sl}(2)}$. If $n=1$ then $fv=0$ and we have
\begin{align*}
  \Cas v = \big(2ef + {\textstyle {1 \over 2}}h(h-2)\big)v = \chi_{\lambda-\rho}(\Cas) v.
\end{align*}
The rest of the proof is by induction on $n$. Let us assume that the statement holds for some $n \in \N$. If for a weight vector $v \in M$ with weight $\lambda \in \mfrak{h}^*$ holds $f^{n+1}v=0$, then $w = fv$ is a weight vector with weight $\lambda-\alpha \in \mfrak{h}^*$ and $f^nw=0$, which by the induction assumption gives us
\begin{align*}
  0 = \prod_{k=1}^n (\Cas - \chi_{\lambda-\rho-k\alpha}(\Cas))w = \prod_{k=1}^n (\Cas - \chi_{\lambda-\rho-k\alpha}(\Cas))fv = f\prod_{k=1}^n (\Cas - \chi_{\lambda-\rho-k\alpha}(\Cas))v.
\end{align*}
Together with the fact ${1 \over 2}h(h-2)v = \chi_{\lambda-\rho}(\Cas) v$, we immediately obtain
\begin{align*}
 0 = \big(2ef + {\textstyle {1 \over 2}}h(h-2) - \chi_{\lambda-\rho}(\Cas)\big) \prod_{k=1}^n (\Cas - \chi_{\lambda-\rho-k\alpha}(\Cas))v = \prod_{k=0}^n (\Cas - \chi_{\lambda-\rho-k\alpha}(\Cas))v.
\end{align*}
This implies the first statement.
\smallskip

ii) Let $v \in M$ be a weight vector with weight $\lambda \in \mfrak{h}^*$. Since $M$ is locally $\widebar{\mfrak{n}}$-finite, there exists $n \in \N$
 such that $f^nv=0$. We have
\begin{align*}
  \prod_{k=0}^{n-1} (z-\chi_{\lambda-\rho-k\alpha}(z))v = 0
\end{align*}
for all $z \in \mfrak{z}_{\mfrak{sl}(2)}$, which implies that $v \in \sum_{k=0}^{n-1} M_{\chi_{\lambda-\rho-k\alpha}}$. By using the fact that $M$ is a weight $\mfrak{sl}(2)$-module we immediately obtain that $M$ is a $\mfrak{z}_{\mfrak{sl}(2)}$-weight module and the $\mfrak{z}_{\mfrak{sl}(2)}$-weights of $M$ may be only of the form $\chi_{\mu-\rho}$ for $\mu \in \mfrak{h}^*$.

Since $\chi_{\mu_1-\rho} = \chi_{\mu_2-\rho}$ for $\mu_1, \mu_2 \in \mfrak{h}^*$ and $\mu_1 \neq \mu_2$ if and only if $\mu_2 -\rho = s(\mu_1-\rho) = -\mu_1 +\rho$, where $s$ is the nontrivial element of the Weyl group $W$ of $\mfrak{sl}(2)$, we get that $\chi_{\lambda-\rho-k_1\alpha} = \chi_{\lambda-\rho-k_2\alpha}$ for $0 \leq k_1 < k_2$ provided $\lambda - \rho =  {1\over 2} (k_1+k_2) \alpha$. Moreover, we have $\lambda - \rho - k_1\alpha = {1\over 2}(k_2-k_1)\alpha$ and $\lambda - \rho - k_2\alpha = -{1\over 2}(k_2-k_1)\alpha$. Therefore, the multiplicity of a $\mfrak{z}_{\mfrak{sl}(2)}$-character $\chi_{\mu-\rho}$ for $\mu \in \mfrak{h}^*$ occurring in the decomposition $\smash{\prod_{k=0}^{n-1}} (z-\chi_{\lambda-\rho-k\alpha}(z))$ is at most $2$ and it happens only if $w(\mu-\rho)-\rho \in \Lambda^+$ for some $w \in W$, otherwise is $1$. This gives us the second statement.}

\vspace{-2mm}


\section{$\alpha$-Gelfand--Tsetlin modules}
\label{sec:Gelfand--Tsetlin}


\subsection{Definitions and examples}

Let us denote by $\mfrak{s}_\alpha$ for a positive root $\alpha \in \Delta^+$ the Lie subalgebra of $\mfrak{g}$ determined through the $\mfrak{sl}(2)$-triple $(e_\alpha,h_\alpha,f_\alpha)$. Then the quadratic Casimir element $\Cas_\alpha$ given by
\begin{align}
   \Cas_\alpha = e_\alpha f_\alpha + f_\alpha e_\alpha + {\textstyle {1 \over 2}} h_\alpha^2
\end{align}
is a free generator of the center $\mfrak{z}_{\mfrak{s}_\alpha}$\! of $U(\mfrak{s}_\alpha)$.
\medskip

\definition{Let us denote by $\Gamma_\alpha$ for $\alpha \in \Delta^+$ the commutative $\C$-subalgebra of $U(\mfrak{g})$  generated by the Cartan subalgebra $\mfrak{h}$ and by the center $\mfrak{z}_{\mfrak{s}_\alpha}$\! of $U(\mfrak{s}_\alpha)$. The $\C$-algebra $\Gamma_\alpha$ is freely generated by the coroots $h_\gamma$ for $\gamma \in \Pi$ and by the quadratic Casimir element $\Cas_\alpha$.}

The $\Gamma_\alpha$-Gelfand--Tsetlin modules  will be simply called \emph{$\alpha$-Gelfand--Tsetlin modules}. When $\alpha$ is the maximal root of a simple Lie algebra $\mfrak{g}$ such modules were considered in \cite{Futorny-Krizka2017}.
\medskip

The following proposition is an immediate consequence of Proposition \ref{prop:Jordan blocks}.
\medskip

\proposition{Let $\alpha \in \Delta^+$ and let $M$ be a weight $\mfrak{g}$-module which is locally $\mfrak{s}_\alpha^-$-finite. Then $M$ is an $\alpha$-Gelfand--Tsetlin module with Jordan blocks of size at most $2$.}

Let us consider a weight $\mfrak{g}$-module $M$. Since $[\mfrak{h}, \mfrak{s}_\alpha^-] = \mfrak{s}_\alpha^-$ for $\alpha \in \Delta^+$, the Lie algebra cohomology groups $H^n(\mfrak{s}_\alpha^-;M)$ are weight $\mfrak{h}$-modules for all $n \in \N_0$. Let us note that $H^n(\mfrak{s}_\alpha^-;M)=0$ for all $n>1$ and $H^0(\mfrak{s}_\alpha^-;M)\simeq \smash{M^{\mfrak{s}_\alpha^-}}$, $H^1(\mfrak{s}_\alpha^-;M) \simeq M/\mfrak{s}_\alpha^-M$.
\medskip

\theorem{\label{thm:GT modules general}
Let $\alpha \in \Delta^+$ and let $M$ be a weight $\mfrak{g}$-module which is locally $\mfrak{s}_\alpha^-$-finite. Then $M$ is an $\alpha$-Gelfand--Tsetlin module with finite $\Gamma_\alpha$-multiplicities if and only if the zeroth cohomology group $H^0(\mfrak{s}_\alpha^-;M)$ is a weight $\mfrak{h}$-module with finite-dimensional weight spaces.}

\proof{As $f_\alpha$ acts locally nilpotently on $M$, by Proposition \ref{prop:Jordan blocks} we obtain that $M$ is a $\mfrak{z}_{\mfrak{s}_\alpha}$-weight module and hence $\Gamma_\alpha$-weight module. Further, since $M_\mu$ for $\mu \in \mfrak{h}^*$ is a $\mfrak{z}_{\mfrak{s}_\alpha}$-module, we only need to show that $M_\mu$ has finite $\mfrak{z}_{\mfrak{s}_\alpha}$-multiplicities if and only if $H^0(\mfrak{s}_\alpha^-;M)$ is a weight $\mfrak{h}$-module with finite-dimensional weight spaces.

Let us define an increasing filtration $\{F_nM_\mu\}_{n \in \N_0}$ of $M_\mu$ by $F_nM_\mu = \{v \in M_\mu;\, f_\alpha^nv=0\}$. Since $F_nM_\mu$ is a $\mfrak{z}_{\mfrak{s}_\alpha}$-module for $n \in \N_0$, we get that $F^n\!M_\mu=F_{n+1}M_\mu/F_nM_\mu$ is a $\mfrak{z}_{\mfrak{s}_\alpha}$-module for $n \in \N_0$. Therefore, for $z \in \mfrak{z}_{\mfrak{s}_\alpha}$ we have the induced linear mapping
\begin{align*}
  \gr^F_n\!z  \colon F^n\!M_\mu \rarr F^n\!M_\mu.
\end{align*}
Moreover, by Proposition \ref{prop:Jordan blocks} we have
\begin{align*}
  (z -\chi_{\mu_\alpha-\rho_\alpha- n\alpha}(z))f_\alpha^nv = 0
\end{align*}
for $v \in F_{n+1}M_\mu$ and $z \in \mfrak{z}_{\mfrak{s}_\alpha}$, where $\rho_\alpha = \rho_{\mfrak{s}_\alpha} ={1 \over 2} \alpha$ and $\mu_\alpha = \mu_{|\mfrak{s}_\alpha \cap \mfrak{h}}$, which immediately implies that $(z-\chi_{\mu_\alpha-\rho_\alpha-n\alpha}(z))v \in F_nM_\mu$. Hence, we have $\gr^F_n\! z_{|F^n\!M_\mu} = \chi_{\mu_\alpha-\rho_\alpha-n\alpha}(z)\, \id_{F^n\!M_\mu}$, which means that $F^n\!M_\mu= (F^n\!M_\mu)_{\chi_{\mu_\alpha-\rho_\alpha-n\alpha}}$. Using the fact that
\begin{align*}
  (F_{n+1}M_\mu)_\chi/(F_nM_\mu)_\chi \simeq (F^n\!M_\mu)_\chi
\end{align*}
for $\chi \in \Hom(\mfrak{z}_{\mfrak{s}_\alpha},\C)$, we obtain
\begin{align*}
  \dim\, (M_\mu)_{\chi_{\mu_\alpha-\rho_\alpha-n\alpha}} =
  \begin{cases}
    \dim F^n\!M_\mu & \text{if $\mu(h_\alpha)-n-1 \notin \N_0 \setminus \{n\}$}, \\
    \dim F^n\!M_\mu + \dim F^{\mu(h_\alpha)-n-1}\!M_\mu & \text{if $\mu(h_\alpha)-n-1 \in \N_0 \setminus \{n\}$}
  \end{cases}
\end{align*}
for $n \in \N_0$, where we used that $\chi_{\mu_\alpha-\rho_\alpha-k\alpha} = \chi_{\mu_\alpha-\rho_\alpha-n\alpha}$ if and only if $k=n$ or $k=\mu(h_\alpha)-n-1$ and $k \in \N_0$.

Now, let us assume that $H^0(\mfrak{s}_\alpha^-;M)\simeq F_1M$ is a weight $\mfrak{h}$-module with finite-dimensional weight spaces. We show by induction on $n$ that $\dim F_nM_\mu < \infty$ for $n \in \N_0$ and $\mu \in \mfrak{h}^*$. For $n=0$ it is trivial and for $n=1$ it follows from the requirement on $H^0(\mfrak{s}_\alpha^-;M)$. Let us assume that it holds for some $n \in \N$ and let us consider the linear mapping $\smash{{f_\alpha^n}_{|F_{n+1}M}} \colon F_{n+1}M \rarr F_1M$. Since $\Ker f_\alpha^n = F_nM$, we obtain that the induced linear mapping $\smash{{f_\alpha^n}_{|F_{n+1}M}} \colon F_{n+1}M/F_nM \rarr F_1M$ is injective. As we have $f_\alpha^n F_{n+1}M_\mu \subset F_1M_{\mu-n\alpha}$, $\dim F_1M_{\mu-n\alpha} < \infty$ and by the induction assumption also $\dim F_nM_\mu < \infty$, we get $\dim F_{n+1}M_\mu < \infty$. Hence, we have $\dim\, (M_\mu)_\chi < \infty$ for $\chi \in \Hom(\mfrak{z}_{\mfrak{s}_\alpha},\C)$.

On the other hand, if $M_\mu$ has finite $\mfrak{z}_{\mfrak{s}_\alpha}$-multiplicities, then we have $\dim\, (F_1M)_\mu = \dim F_1M_\mu = \dim F^1\!M_\mu \leq \dim\, (M_\mu)_{\chi_{\mu_\alpha-\rho_\alpha-\alpha}} < \infty$. Hence, we have that  $H^0(\mfrak{s}_\alpha^-;M) \simeq F_1M$ is a weight $\mfrak{h}$-module with finite-dimensional weight spaces.}

By duality we can analogously prove a similar statement if we replace $f_\alpha$ by $e_\alpha$ and $H^0(\mfrak{s}_\alpha^-;M)$ by $H^0(\mfrak{s}_\alpha^+;M)$.


\subsection{Categories of $\alpha$-Gelfand--Tsetlin modules}

We denote by $\mcal{M}(\mfrak{g})$ the category of $\mfrak{g}$-modules and by $\mcal{E}(\mfrak{g})$ the category of finite-dimensional $\mfrak{g}$-modules. For a Lie subalgebra $\mfrak{a}$ of $\mfrak{g}$ we introduce the full subcategory $\mcal{I}(\mfrak{g},\mfrak{a})$ of
$\mcal{M}(\mfrak{g})$ consisting of locally $\mfrak{a}$-finite weight $\mfrak{g}$-modules and the full subcategory $\mcal{I}_f(\mfrak{g},\mfrak{a})$ of $\mcal{M}(\mfrak{g})$ of finitely generated locally $\mfrak{a}$-finite weight $\mfrak{g}$-modules. Therefore, for $\alpha \in \Delta^+$ we have the full subcategories $\mcal{I}(\mfrak{g},\mfrak{s}_\alpha^+)$ and $\mcal{I}(\mfrak{g},\mfrak{s}_\alpha^-)$ of $\mcal{M}(\mfrak{g})$ assigned to the Lie subalgebras $\mfrak{s}_\alpha^+ = \mfrak{s}_\alpha \cap \mfrak{n}$ and $\mfrak{s}_\alpha^- = \mfrak{s}_\alpha \cap \widebar{\mfrak{n}}$ of $\mfrak{g}$.
\medskip

The following statement is obvious, see Proposition 1.3 in \cite{Futorny-Krizka2017} for details.
\medskip

\proposition{\label{prop-subcat}
 \begin{enumerate}[topsep=3pt,itemsep=0pt]
\item[i)] The categories $\mcal{I}(\mfrak{g},\mfrak{s}_\alpha^+)$ and $\mcal{I}(\mfrak{g},\mfrak{s}_\alpha^-)$ are full subcategories of $\mcal{H}(\mathfrak{g}, \Gamma_\alpha)$ for $\alpha \in \Delta^+$.
\item[ii)] Let $M \in \mcal{M}(\mfrak{g})$ be a simple weight $\mfrak{g}$-module and let $\alpha \in \Delta^+$. If there exists a positive integer $n \in \N$ and a nonzero vector $v \in M$ such that $e_\alpha^n v=0$ (or $f_\alpha^n v=0$), then  $M \in \mcal{I}(\mfrak{g},\mfrak{s}_\alpha^+)$ (or $M \in  \mcal{I}(\mfrak{g},\mfrak{s}_\alpha^-)$).
 \end{enumerate}}

We introduce $\mcal{I}_{\rm fin}(\mfrak{g},\mfrak{s}_\alpha^+) = \mcal{I}_f(\mfrak{g},\mfrak{s}_\alpha^+) \cap \mcal{H}_{\rm fin}(\mathfrak{g}, \Gamma_\alpha)$ and $\mcal{I}_{\rm fin}(\mfrak{g},\mfrak{s}_\alpha^-) = \mcal{I}_f(\mfrak{g},\mfrak{s}_\alpha^-) \cap \mcal{H}_{\rm fin}(\mathfrak{g}, \Gamma_\alpha)$ for $\alpha \in \Delta^+$.
\medskip

\remark{
\begin{itemize}[topsep=3pt,itemsep=0pt]
\item We have $\mcal{E}(\mfrak{g})\subset \mcal{I}_{\rm fin}(\mfrak{g},\mfrak{s}_\alpha^+) \cap \mcal{I}_{\rm fin}(\mfrak{g},\mfrak{s}_\alpha^-)$. Moreover, it holds $\mcal{E}(\mfrak{g})= \mcal{I}_{\rm fin}(\mfrak{g},\mfrak{s}_\alpha^+) \cap \mcal{I}_{\rm fin}(\mfrak{g},\mfrak{s}_\alpha^-)$ provided $\mathfrak{g}=\mfrak{sl}(2)$;
 \item The category $\mcal{I}_{\rm fin}(\mfrak{g},\mfrak{s}_\alpha^+)$ contains the category $\mcal{O}(\mfrak{g})$.
 \item Examples of simple $\mfrak{g}$-modules with infinite-dimensional weight spaces but with finite $\Gamma_\alpha$-multiplicities were constructed in \cite{Futorny-Krizka2017} for the maximal root  $\alpha$  of a simple Lie algebra $\mfrak{g}$.
\end{itemize}}

\proposition{\label{prop:tensor category}
Let $\mfrak{a}$ be a nilpotent Lie subalgebra of $\mfrak{g}$. Then the category $\mcal{I}(\mfrak{g},\mfrak{a})$ is a tensor category.}

\proof{Since $\mfrak{a}$ is a nilpotent Lie algebra, the condition that a $\mfrak{g}$-module $M$ is locally $\mfrak{a}$-finite is equivalent to saying that each element $a \in \mfrak{a}$ acts locally nilpotently on $M$. The rest of the proof follows immediately from the definition of the tensor product of $\mfrak{g}$-modules.}

Let us note that the categories $\mcal{H}(\mfrak{g}, \Gamma_\alpha)$, $\mcal{H}_{\rm fin}(\mfrak{g}, \Gamma_\alpha)$, $\mcal{I}_f(\mfrak{g},\mfrak{s}_\alpha^+)$, $\mcal{I}_f(\mfrak{g},\mfrak{s}_\alpha^-)$, $\mcal{I}_{\rm fin}(\mfrak{g},\mfrak{s}_\alpha^+)$, $\mcal{I}_{\rm fin}(\mfrak{g},\mfrak{s}_\alpha^-)$ and $\mcal{I}_{\rm fin}(\mfrak{g},\mfrak{s}_\alpha^+) \cap \mcal{I}_{\rm fin}(\mfrak{g},\mfrak{s}_\alpha^-)$ are not tensor categories (they are not closed with respect to the tensor product of $\mfrak{g}$-modules) for $\alpha \in \Delta^+$, except the case $\mathfrak{g}=\mfrak{sl}(2)$ when $\mcal{I}_{\rm fin}(\mfrak{g},\mfrak{s}_\alpha^+) \cap \mcal{I}_{\rm fin}(\mfrak{g},\mfrak{s}_\alpha^-) = \mcal{E}(\mfrak{g})$. On the other hand, by Proposition \ref{prop:tensor category} we have that $\mcal{I}(\mfrak{g},\mfrak{s}_\alpha^+)$ and $\mcal{I}(\mfrak{g},\mfrak{s}_\alpha^-)$ are tensor categories.
\medskip

Let $\Gamma$ be a commutative $\C$-subalgebra of $U(\mfrak{g})$ containing $\Gamma_\alpha$. Then we have the following chain of embeddings
\begin{align*}
 \mcal{H}_{\rm fin}(\mfrak{g},\Gamma_\alpha) \subset \mcal{H}_{\rm fin}(\mfrak{g},\Gamma) \subset \mcal{H}(\mfrak{g},\Gamma) \subset \mcal{H}(\mfrak{g},\Gamma_\alpha)
\end{align*}
of categories. In particular, we see that $\mcal{I}_{\rm fin}(\mfrak{g},\mfrak{s}_\alpha^\pm)$ is a full subcategory of $\mcal{H}_{\rm fin}(\mathfrak{g}, \Gamma)$. However, ${\mcal{I}}(\mfrak{g},\mfrak{s}_\alpha^{\pm})$ and $\mcal{I}_f(\mfrak{g},\mfrak{s}_\alpha^{\pm})$ are not subcategories of $\mcal{H}(\mathfrak{g}, \Gamma)$ in general. Also let us note that neither $\mcal{H}(\mfrak{g}, \Gamma)$ nor its subcategories $\mcal{H}(\mfrak{g}, \Gamma)\cap \mcal{I}(\mfrak{g},\mfrak{s}_\alpha^+)$ and $\mcal{H}(\mfrak{g}, \Gamma)\cap \mcal{I}(\mfrak{g},\mfrak{s}_\alpha^-)$ are tensor categories if $\mfrak{g}\neq \mfrak{sl}(2)$. Nevertheless, we have the following result for $ \mfrak{sl}(3)$.
\medskip

\theorem{\label{thm-sl3}
Let $\mfrak{g} = \mfrak{sl}(3)$ and $\alpha \in \Delta^+$. Let $V$ and $W$ be simple $\mfrak{g}$-modules in $ \mcal{I}(\mfrak{g},\mfrak{s}_\alpha^+)$ (or in $\mcal{I}(\mfrak{g},\mfrak{s}_\alpha^-)$). Then the modules $V$, $W$ and every simple subquotient of $V\otimes_\C W$ belong to the category $\mcal{I}_{\rm fin}(\mfrak{g},\mfrak{s}_\alpha^+)$ (or to the category $\mcal{I}_{\rm fin}(\mfrak{g},\mfrak{s}_\alpha^-)$) and hence to the category  $\mcal{H}_{\rm fin}(\mfrak{g}, \Gamma)$ for any commutative $\C$-subalgebra $\Gamma$ of $U(\mfrak{g})$ containing $\Gamma_\alpha$.}

\proof{Indeed, as $\Gamma$ is a commutative $\C$-subalgebra containing $\Gamma_{\alpha}$, then either $\Gamma=\Gamma_\alpha$ or it is generated over $\Gamma_\alpha$  by central elements of $U(\mfrak{g})$. Then the simplicity of modules and Proposition \ref{prop-subcat} (i)
imply that $V$ and $W$ belong to the category $\mcal{H}(\mfrak{g}, \Gamma)$ for any commutative $\C$-subalgebra $\Gamma$ of $U(\mfrak{g})$ containing $\Gamma_\alpha$. In particular, we have $V, W\in \mcal{H}(\mfrak{g}, \Gamma)$, where $\Gamma$ is a Gelfand--Tsetlin subalgebra generated by
$\Gamma_{\alpha} $ and by the center of $U(\mfrak{g})$. The structure of simple $\Gamma$-Gelfand--Tsetlin modules for $\mfrak{sl}(3)$ was described in \cite{FGR-sl(3)}, see also \cite{Futorny-Grantcharov-Ramirez2014}. It follows that $V, W\in \mcal{I}_{\rm fin}(\mfrak{g},\mfrak{s}_\alpha^+)$ (or that $V, W\in \mcal{I}_{\rm fin}(\mfrak{g},\mfrak{s}_\alpha^-)$). Every simple subquotient $U$ of $V\otimes_\C W$ is clearly a weight module with a locally finite action of $\mfrak{s}_\alpha^+$ (or $\mfrak{s}_\alpha^-$). Again, the results of \cite{FGR-sl(3)} imply that  $\Gamma_{\alpha}$-multiplicities in $U$ are finite, and the statement follows.}

Theorem \ref{thm-sl3} suggests a possible generalization for $\mfrak{sl}(n)$. Indeed, let us consider the simple Lie algebra $\mfrak{g}=\mfrak{sl}(n)$ with a triangular decomposition $\mfrak{g} = \widebar{\mfrak{n}} \oplus \mfrak{h} \oplus \mfrak{n}$ and with the set of simple roots $\Pi = \{\alpha_1,\alpha_2,\dots,\alpha_{n-1}\}$ such that the corresponding Cartan matrix $A=(a_{ij})_{1\leq i,j \leq n-1}$ is given by $a_{ii} =2$, $a_{ij}= -1$ if $|i-j|=1$ and $a_{ij}=0$ if $|i-j| \geq 2$. Further, let us denote by $\mfrak{g}_k$ for $k=1,2,\dots,n$ the Lie subalgebra of $\mfrak{g}$ generated by the root subspaces $\mfrak{g}_{\alpha_1}, \dots, \mfrak{g}_{\alpha_{k-1}}$ and $\mfrak{g}_{-\alpha_1}, \dots, \mfrak{g}_{-\alpha_{k-1}}$. Then we obtain a finite sequence
\begin{align*}
  0 = \mfrak{g}_1 \subset \mfrak{g}_2 \subset \dots \subset \mfrak{g}_{n-1} = \mfrak{g}
\end{align*}
of Lie subalgebras of $\mfrak{g}$ such that $\mfrak{g}_k \simeq \mfrak{sl}(k)$ for $k=2,3,\dots,n$. We have also the induced triangular decomposition
\begin{align*}
  \mfrak{g}_k = \widebar{\mfrak{n}}_k \oplus \mfrak{h}_k \oplus \mfrak{n}_k
\end{align*}
of the Lie algebra $\mfrak{g}_k$ for $k=2,3,\dots,n$, where $\widebar{\mfrak{n}}_k = \widebar{\mfrak{n}} \cap \mfrak{g}_k$, $\mfrak{n}_k = \mfrak{n} \cap \mfrak{g}_k$ and $\mfrak{h}_k = \mfrak{h} \cap \mfrak{g}_k$. Besides, we have a sequence $U(\mfrak{g}_2) \subset U(\mfrak{g}_3) \subset \dots \subset U(\mfrak{g}_n)$ of $\C$-subalgebras of the universal enveloping algebra $U(\mfrak{g})$. Let us denote by $\mfrak{z}_{\mfrak{g}_k}$\! the center of $U(\mfrak{g}_k)$ for $k=2,3,\dots,n$. Then the Gelfand--Tsetlin subalgebra $\Gamma$ of $U(\mfrak{g})$ is generated by $\mfrak{z}_{\mfrak{g}_k}$ for $k=2,3,\dots,n$ and by the Cartan subalgebra $\mfrak{h}$, cf.\ \cite{Drozd-Futorny-Ovsienko1994}. It is a maximal commutative $\C$-subalgebra of $U(\mfrak{g})$.

The following theorem provides a  generalization of Theorem \ref{thm-sl3} to the Lie algebra $\mfrak{sl}(n)$.
\medskip

\theorem{\label{thm-sln}
Let $\mfrak{g} = \mfrak{sl}(n)$.
  \begin{enumerate}[topsep=3pt,itemsep=0pt]
  \item[i)] For $k=2,3,\dots, n$, the categories $\mcal{I}(\mfrak{g},\mfrak{n}_k)$ and $\mcal{I}(\mfrak{g},\widebar{\mfrak{n}}_k)$ are tensor categories.
  \item[ii)] Let $V$ and $W$ be simple $\mathfrak{g}$-modules in  $\mcal{I}(\mfrak{g},\mfrak{n}_{n-1})$ or in  $\mcal{I}(\mfrak{g},\widebar{\mfrak{n}}_{n-1})$. Then the $\mfrak{g}$-modules $V$, $W$ and every simple subquotient of $V\otimes_\C W$ belongs to $\mcal{H}(\mathfrak{g}, \Gamma)$.
  \end{enumerate}
}

\proof{The proof is analogous to the proof of Theorem \ref{thm-sl3}.}

\vspace{-2mm}


\section{Twisting functors and generalized Verma modules}
\label{sec:Twisting functor and Verma modules}

We define the twisting functor $T_\alpha$ assigned to a positive root $\alpha \in \Delta^+$ of a semisimple Lie algebra $\mfrak{g}$ and describe main properties of this functor. By applying of the twisting functor $T_\alpha$ to generalized Verma modules we construct $\alpha$-Gelfand--Tsetlin $\mfrak{g}$-modules with finite $\Gamma_\alpha$-multiplicities. If $\alpha$ is a simple root, then these modules are twisted Verma modules up to conjugation of the action of $\mfrak{g}$, see \cite{Andersen-Lauritzen2003}, \cite{Krizka-Somberg2015b}, \cite{Musson2019}.


\subsection{Twisting functor $T_\alpha$}
\label{subsec:Twisting functor}

We use the notation introduced in  previous sections. Let $e_\alpha \in \mfrak{g}_\alpha$ and $f_\alpha \in \mfrak{g}_{-\alpha}$ for $\alpha \in \Delta^+$ be nonzero elements of $\mfrak{g}$ satisfying $[e_\alpha,f_\alpha]=h_\alpha$. The multiplicative set $\{f_\alpha^n;\, n \in \N_0\}$ in the universal enveloping algebra $U(\mfrak{g})$ of $\mfrak{g}$ is a left (right) denominator set, since $f_\alpha$ is a locally $\ad$-nilpotent regular element. Therefore, based on the previous general construction we have defined the \emph{twisting functor}
\begin{align*}
  T_\alpha=T_{f_\alpha} \colon \mcal{M}(\mfrak{g}) \rarr \mcal{M}(\mfrak{g}).
\end{align*}
Besides, we introduce the \emph{partial Zuckerman functor}
\begin{align*}
  S_\alpha \colon \mcal{M}(\mfrak{g}) \rarr \mcal{M}(\mfrak{g})
\end{align*}
defined by
\begin{align*}
  S_\alpha(M) = \{v \in M;\, (\exists n \in \N)\, f_\alpha^n v =0\}
\end{align*}
for $M \in \mcal{M}(\mfrak{g})$.
\medskip

In the next we prove some basic characteristics of the twisting functor $T_\alpha$ for $\alpha \in \Delta^+$. Let us note that the twisting functor $T_\alpha$ is usually considered only for a simple root $\alpha \in \Pi$ in the literature (see e.g.\ \cite{Andersen-Stroppel2003}), which is caused by the fact that $T_\alpha$ preserves the category $\mcal{O}(\mfrak{g})$ up to conjugation of the action of $\mfrak{g}$.
\medskip

For $\alpha \in \Delta^+$ let us consider the subset
\begin{align*}
  \Phi_\alpha = \{\gamma \in \Delta \setminus \{\pm \alpha\};\, \mfrak{g}_{\gamma,\alpha} \subset \mfrak{n}\}
\end{align*}
of $\Delta$, where $\mfrak{g}_{\gamma,\alpha}$ for $\gamma \in \Delta$ is the simple finite-dimensional $\mfrak{s}_\alpha$-module given by
\begin{align*}
  \mfrak{g}_{\gamma,\alpha} = \bigoplus_{j \in \Z} \mfrak{g}_{\gamma+j\alpha}.
\end{align*}
It easily follows that $\Phi_\alpha$ is closed, i.e.\ if $\gamma_1, \gamma_2 \in \Phi_\alpha$ and $\gamma_1+\gamma_2 \in \Delta$, then we have $\gamma_1+\gamma_2 \in \Phi_\alpha$, and that $\Phi_\alpha \subset \Delta^+$. Therefore, the subset $\Phi_\alpha$ of $\Delta^+$ gives rise to the Lie subalgebras $\mfrak{t}_\alpha^+$, $\mfrak{t}_\alpha^-$ and $\mfrak{t}_\alpha$ of $\mfrak{g}$ defined by
\begin{align*}
  \mfrak{t}_\alpha^+ = \mfrak{g}_\alpha \oplus \bigoplus_{\gamma \in \Phi_\alpha} \mfrak{g}_\gamma, \qquad \mfrak{t}_\alpha = \bigoplus_{\gamma \in \Phi_\alpha} \mfrak{g}_\gamma, \qquad  \mfrak{t}_\alpha^- = \mfrak{g}_{-\alpha} \oplus \bigoplus_{\gamma \in \Phi_\alpha} \mfrak{g}_\gamma.
\end{align*}
As we have $\mfrak{s}_\alpha^\pm \subset \mfrak{t}_\alpha^\pm$ for $\alpha \in \Delta^+$, we get immediately  $\mcal{I}(\mfrak{g},\mfrak{t}_\alpha^\pm) \subset \mcal{I}(\mfrak{g},\mfrak{s}_\alpha^\pm)$. Besides, we have
\begin{align*}
\Ad(\dot{s}_\alpha)(\mfrak{s}_\alpha^\pm) = \mfrak{s}_\alpha^\mp \qquad \text{and} \qquad \Ad(\dot{s}_\alpha)(\mfrak{t}_\alpha^\pm) = \mfrak{t}_\alpha^\mp
\end{align*}
for $\alpha \in \Delta^+$, where $\dot{s}_\alpha \in N_G(H)$ is a representative of the element $s_\alpha \in W \simeq N_G(H)/H$, where $G$ and $H$ are connected algebraic groups with Lie algebras $\mfrak{g}$ and $\mfrak{h}$, respectively. The inclusions $\mfrak{s}_\alpha^+ \subset \mfrak{t}_\alpha^+ \subset \mfrak{n}$ of Lie algebras give rise to the embeddings
\begin{align}
  \mcal{O}(\mfrak{g}) \subset \mcal{I}_f(\mfrak{g},\mfrak{t}_\alpha^+) \subset \mcal{I}_f(\mfrak{g},\mfrak{s}_\alpha^+)
\end{align}
of categories for $\alpha \in \Delta^+$.
\medskip

\proposition{\label{prop:locally nilpotent}
Let $M$ be a $\mfrak{g}$-module. Then $T_\alpha(M)$ is a locally $\mfrak{s}_\alpha^-$-finite $\mfrak{g}$-module for any $\alpha \in \Delta^+$. Moreover, if $M$ is a locally $\mfrak{s}_\alpha^-$-finite $\mfrak{g}$-module, then  $T_\alpha(M)=0$.}

\proof{Let $M$ be a $\mfrak{g}$-module. Then by definition we have
\begin{align*}
  T_\alpha(M) = (U(\mfrak{g})_{(f_\alpha)}/U(\mfrak{g})) \otimes_{U(\mfrak{g})} M \simeq M_{(f_\alpha)}/M
\end{align*}
for $\alpha \in \Delta^+$. Since every element of $M_{(f_\alpha)}$ can be written in the form $f_\alpha^{-n} v$ for $n \in \N_0$ and $v \in M$, we obtain immediately that $T_\alpha(M)$ is locally $\mfrak{s}_\alpha^-$-finite. Further, let us assume that $M$ is a locally $\mfrak{s}_\alpha^-$-finite $\mfrak{g}$-module. Since for each $v \in M$ there exists $n_v \in \N_0$ such that $f_\alpha^{n_v}v=0$, we may write $f_\alpha^{-n}v = f_\alpha^{-n-n_v}f_\alpha^{n_v}v = 0$ for $n\in \N_0$. This implies the required statement.}

\lemma{\label{lem:sl2 action}
Let $\alpha \in \Delta^+$. Then we have
\begin{align} \label{eq:sl2 action}
\begin{aligned}
  e_\alpha f_\alpha^{-n} &= f_\alpha^{-n} e_\alpha -nf_\alpha^{-n-1} h_\alpha - n(n+1)f_\alpha^{-n-1}, \\
  h f_\alpha^{-n} &= f_\alpha^{-n} h + n\alpha(h) f_\alpha^{-n}, \\
  f_\alpha f_\alpha^{-n} &= f_\alpha^{-n} f_\alpha
\end{aligned}
\end{align}
in $U(\mfrak{g})_{(f_\alpha)}$ for $n \in \Z$ and $h \in \mfrak{h}^*$.}

\proof{It follows immediately from the formula
\begin{align*}
  a f_\alpha^{-n} = f_\alpha^{-n} \sum_{k=0}^\infty \binom{n+k-1}{k} f_\alpha^{-k} \ad(f_\alpha)^k(a)
\end{align*}
in $U(\mfrak{g})_{(f_\alpha)}$ for all $a \in U(\mfrak{g})$ and $n \in \N_0$.}

\theorem{\label{thm:finitely generated}
We have
\begin{enumerate}[topsep=3pt,itemsep=0pt]
  \item[i)] if $M \in \mcal{I}_f(\mfrak{g},\mfrak{s}_\alpha^+)$, then $T_\alpha(M) \in \mcal{I}_f(\mfrak{g},\mfrak{s}_\alpha^-)$;
  \item[ii)] if $M \in \mcal{I}_f(\mfrak{g},\mfrak{t}_\alpha^+)$, then $T_\alpha(M) \in \mcal{I}_f(\mfrak{g},\mfrak{t}_\alpha^-)$
\end{enumerate}
for $\alpha \in \Delta^+$. Therefore, we have the restricted functors
\begin{align*}
  T_\alpha \colon \mcal{I}_f(\mfrak{g},\mfrak{s}_\alpha^+) \rarr \mcal{I}_f(\mfrak{g},\mfrak{s}_\alpha^-) \qquad \text{and} \qquad T_\alpha \colon \mcal{I}_f(\mfrak{g},\mfrak{t}_\alpha^+) \rarr \mcal{I}_f(\mfrak{g},\mfrak{t}_\alpha^-)
\end{align*}
for $\alpha \in \Delta^+$.}

\proof{i) If $M$ is a weight $\mfrak{g}$-module, then it easily follows from definition that $T_\alpha(M)$ is also a weight $\mfrak{g}$-module. Moreover, by Proposition \ref{prop:locally nilpotent} we have that $T_\alpha(M)$ is locally $\mfrak{s}_\alpha^-$-finite. Hence, the rest of the proof is to show that $T_\alpha(M)$ is finitely generated provided $M$ is finitely generated and locally $\mfrak{s}_\alpha^+$-finite.

Let $R \subset M$ be a finite set of generators of $M$. Then the vector subspace $V = U(\mfrak{s}_\alpha^+)\langle R \rangle$ of $M$ is finite dimensional. Further, let us introduce a filtration $\{F_k V\}_{k \in \N_0}$ on $V$ by
\begin{align*}
  F_k V = \{v \in V;\, e_\alpha^kv=0\}
\end{align*}
for $k \in \N_0$ and let $n_0 \in \N$ be the smallest positive integer satisfying
\begin{align*}
  n_0 \geq \max\{-\mu(h_\alpha);\, \mu \in \mfrak{h}^*,\, V_\mu \neq \{0\},\, \mu(h_\alpha) \in \R\}.
\end{align*}
Let us consider a vector $v \in F_kV \cap V_\mu$ for $\mu \in \mfrak{h}^*$ and $k \in \N_0$. Then by Lemma \ref{lem:sl2 action} we obtain
\begin{align*}
   e_\alpha f_\alpha^{-(n_0+n)}v = f_\alpha^{-(n_0+n)}e_\alpha v - (n_0+n)(\mu(h_\alpha)+n_0+n+1) f_\alpha^{-(n_0+n+1)}v
\end{align*}
for $n \in \N_0$, which together with $f_\alpha f_\alpha^{-n-1}v = f_\alpha^{-n}v$ for $n\in \N_0$ gives us
\begin{align*}
  U(\mfrak{s}_\alpha)f_\alpha^{-n_0}F_kV/\C[f_\alpha^{-1}]F_{k-1}V = \C[f_\alpha^{-1}]F_kV/\C[f_\alpha^{-1}]F_{k-1}V
\end{align*}
for all $k \in \N$. As $F_0V = \{0\}$ and $V$ is a finite-dimensional vector space, we immediately get $U(\mfrak{s}_\alpha)f_\alpha^{-n_0}V = \C[f_\alpha^{-1}]V$.

Further, since for any $a \in U(\mfrak{g})$ and $n \in \N_0$ there exist $b \in U(\mfrak{g})$ and $m \in \N_0$ satisfying $f_\alpha^{-n}a = bf_\alpha^{-m}$, we obtain $f_\alpha^{-n}U(\mfrak{g}) \subset U(\mfrak{g})\C[f_\alpha^{-1}] \subset U(\mfrak{g})_{(f_\alpha)}$ for $n \in \N_0$. Hence, we may write
\begin{align*}
  f_\alpha^{-n}M = f_\alpha^{-n}U(\mfrak{g})V \subset U(\mfrak{g})\C[f_\alpha^{-1}]V,
\end{align*}
which implies $M_{(f_\alpha)} =U(\mfrak{g})\C[f_\alpha^{-1}]V = U(\mfrak{g})f_\alpha^{-n_0}V$. In other words, this means that $T_\alpha(M) \simeq M_{(f_\alpha)}/M$ is finitely generated and the number of generators is bounded by $\dim V$.
\smallskip

ii) Since $\mcal{I}_f(\mfrak{g},\mfrak{t}_\alpha^+)$ is a full subcategory of $\mcal{I}_f(\mfrak{g},\mfrak{s}_\alpha^+)$, by item (i) we only need to show that $T_\alpha(M)$ is locally $\mfrak{t}_\alpha^-$-finite if $M$ is locally $\mfrak{t}_\alpha^+$-finite. As $\mfrak{t}_\alpha^-$ is a nilpotent Lie algebra and $T_\alpha(M)$ is locally $\mfrak{s}_\alpha^-$-finite the condition that $T_\alpha(M)$ is locally $\mfrak{t}_\alpha^-$-finite is equivalent to saying that the element $e_\gamma$, where $\mfrak{g}_{-\gamma} = \C e_\gamma$, acts locally nilpotently on $T_\alpha(M)$ for $\gamma \in \Phi_\alpha$.

Let $\htt \colon \Z\Delta \rarr \Z$ be the $\Z$-linear height function with $\htt(\alpha)=1$ if $\alpha$ is a simple root. Then we get an $\N_0$-grading $U(\mfrak{t}_\alpha) = \bigoplus_{n \in \N_0} U(\mfrak{t}_\alpha)_n$, where
\begin{align*}
  U(\mfrak{t}_\alpha)_n = \bigoplus_{\mu \in \mfrak{h}^*\!,\, \htt(\mu)=n} U(\mfrak{t}_\alpha)_\mu.
\end{align*}
Let $\gamma \in \Phi_\alpha$ and let $r \in \N_0$ be the smallest nonnegative integer such that $\gamma-(r+1)\alpha \notin \Phi_\alpha$. Let us recall that if $\gamma - k\alpha \in \Delta$ for some $k \in \Z$ then $\gamma - k\alpha \in \Phi_\alpha$. Let us consider a vector $v \in M$. Then from the formula
\begin{align*}
  a f_\alpha^{-n} = f_\alpha^{-n} \sum_{k=0}^\infty \binom{n+k-1}{k} f_\alpha^{-k} \ad(f_\alpha)^k(a)
\end{align*}
in $U(\mfrak{g})_{(f_\alpha)}$ for $a \in U(\mfrak{g})$ and $n \in \N_0$, we obtain
\begin{align*}
  e_\gamma^t f_\alpha^{-n}v = f_\alpha^{-n} \sum_{k=0}^{tr} \binom{n+k-1}{k} f_\alpha^{-k} \ad(f_\alpha)^k(e_\gamma^t)v
\end{align*}
for $t \in \N$ and $n \in \N_0$, where we used the fact that $\ad(f_\alpha)^k(e_\gamma) \neq 0$ only for $k=0,1,\dots,r$. Moreover, we have $\ad(f_\alpha)^k(e_\gamma^t) \in U(\mfrak{t}_\alpha)_{t\gamma-k\alpha} \subset U(\mfrak{t}_\alpha)_{\htt(t\gamma-k\alpha)}$ for $k=0,1,\dots,tr$. Since $M$ is locally $\mfrak{t}_\alpha^+$-finite, there exists an integer $n_v \in \N_0$ such that $U(\mfrak{t}_\alpha)_nv=\{0\}$ for $n > n_v$. Therefore, it is enough to show that $\htt(t\gamma-k\alpha) > n_v$ for $k=0,1,\dots,tr$.

As we may write
\begin{align*}
  \htt(t\gamma-k\alpha) \geq \htt(t\gamma-tr\alpha) = t\htt(\gamma-r\alpha) \geq t
\end{align*}
for $k = 0,1,\dots,tr$, since $\gamma - r \alpha \in \Phi_\alpha$ and hence $\htt(\gamma -r\alpha) \geq 1$, we obtain that $e_\gamma^t f_\alpha^{-n}v=0$ for $n \in \N_0$ provided $t > n_v$. Hence, the element $e_\gamma$ acts locally nilpotently on $T_\alpha(M)$ for $\gamma \in \Phi_\alpha$.}

The following shows how the twisting functor $T_{\alpha}$ can be used to construct $\alpha$-Gelfand--Tsetlin modules with finite $\Gamma_\alpha$-multiplicities.
\medskip

\theorem{\label{thm:highest weight modules condition}
Let $M$ be a weight $\mfrak{g}$-module and let $\alpha \in \Delta^+$.
\begin{itemize}[topsep=3pt,itemsep=0pt]
\item[i)] The $\mfrak{g}$-module $T_\alpha(M)$ is an $\alpha$-Gelfand--Tsetlin module with finite $\Gamma_\alpha$-multiplicities if and only if the first cohomology group $H^1(\mfrak{s}_\alpha^-;M)$ is a weight $\mfrak{h}$-module with finite-dimensional weight spaces.
\item[ii)] If $M$ is a highest weight $\mfrak{g}$-module then $T_\alpha(M)$ is a locally $\mfrak{t}_\alpha^-$-finite cyclic weight $\mfrak{g}$-module with finite $\Gamma_\alpha$-multiplicities.
\end{itemize}
}

\proof{i) Let $M$ be a weight $\mfrak{g}$-module. Then by Proposition \ref{prop:locally nilpotent} we have that $T_\alpha(M)$ is a locally $\mfrak{s}_\alpha^-$-finite weight $\mfrak{g}$-module for $\alpha \in \Delta^+$. Hence, we may apply Theorem \ref{thm:GT modules general} on $T_\alpha(M)$ and we obtain that $T_\alpha(M)$ is an $\alpha$-Gelfand--Tsetlin module with finite $\Gamma_\alpha$-multiplicities if and only if $H^0(s_\alpha^-;T_\alpha(M))$ is a weight $\mfrak{h}$-module with finite-dimensional weight spaces. Further, the $\C$-linear mapping
$\varphi_\alpha \colon M \rarr T_\alpha(M)$ defined by
\begin{align*}
  \varphi_\alpha(v) = f_\alpha^{-1}v
\end{align*}
for $v \in M$ gives rise to the $\C$-linear mapping
\begin{align*}
 \widetilde{\varphi}_\alpha \colon H^1(\mfrak{s}_\alpha^-;M) \rarr H^0(\mfrak{s}_\alpha^-;T_\alpha(M))
\end{align*}
for $\alpha \in \Delta^+$, which is in fact an isomorphism. Therefore, we obtain an isomorphism $H^1(\mfrak{s}_\alpha^-;M) \simeq H^0(\mfrak{s}_\alpha^-;T_\alpha(M)) \otimes_\C \C_{-\alpha}$ of $\mfrak{h}$-modules, where $\C_{-\alpha}$ is the $1$-dimensional $\mfrak{h}$-module determined by the character $-\alpha$ of $\mfrak{h}$, which implies the first statement.
\smallskip

ii) If $M$ is a highest weight $\mfrak{g}$-module, then it belongs to the category $\mcal{O}(\mfrak{g})$. Hence, using Theorem \ref{thm:finitely generated} we obtain that $T_\alpha(M)$ is a locally $\mfrak{t}_\alpha^-$-finite weight $\mfrak{g}$-module. In fact, from the proof of Theorem \ref{thm:finitely generated} (i) it follows that $T_\alpha(M)$ is not only finitely generated by also cyclic. To finish the proof, we need to show by Theorem \ref{thm:highest weight modules condition} that $H^1(\mfrak{s}_\alpha^-;M)$ is a weight $\mfrak{h}$-module with finite-dimensional weight spaces. Since $H^1(\mfrak{s}_\alpha^-;M) \simeq M/\mfrak{s}_\alpha^-M$ and $M$ is a weight $\mfrak{h}$-module with finite-dimensional weight spaces, we immediately obtain that also $M/\mfrak{s}_\alpha^-M$ is a weight $\mfrak{h}$-module with finite-dimensional weight spaces. This gives us the required statement.}

If we denote by $\Theta_\alpha \colon \mcal{M}(\mfrak{g}) \rarr \mcal{M}(\mfrak{g})$ the functor sending $\mfrak{g}$-module to the same $\mfrak{g}$-module with the action twisted by the automorphism $\Ad(\dot{s}_\alpha) \colon \mfrak{g} \rarr \mfrak{g}$, then we obtain the endofunctor
\begin{align*}
 \Theta_\alpha \circ T_\alpha \colon \mcal{I}_f(\mfrak{g},\mfrak{t}_\alpha^+) \rarr \mcal{I}_f(\mfrak{g},\mfrak{t}_\alpha^+).
\end{align*}
Moreover, for $\alpha \in \Pi$ we have $\mfrak{t}_\alpha^+ = \mfrak{n}$ which implies $\mcal{I}_f(\mfrak{g},\mfrak{t}_\alpha^+) = \mcal{O}(\mfrak{g})$ and in this case the functor coincides with the Arkhipov's twisting functor, see \cite{Arkhipov1997}, \cite{Arkhipov2004}. On the other hand, we have $\mcal{I}_f(\mfrak{g},\mfrak{t}_\alpha^+) = \mcal{I}_f(\mfrak{g},\mfrak{s}_\alpha^+)$ if $\alpha$ is the maximal root of a simple Lie algebra $\mfrak{g}$.


\subsection{Left derived functor of $T_\alpha$}

Let us recall that we have the following obvious inclusions
\begin{align*}
  \mcal{O}(\mfrak{g}) \subset \mcal{I}_f(\mfrak{g},\mfrak{t}_\alpha^+) \subset \mcal{I}_f(\mfrak{g},\mfrak{s}_\alpha^+)
\end{align*}
of categories for $\alpha \in \Delta^+$. For that reason, we may restrict the twisting functor $T_\alpha$ for $\alpha \in \Delta^+$ on the category $\mcal{O}(\mfrak{g})$. Since $T_\alpha \colon \mcal{O}(\mfrak{g}) \rarr \mcal{I}_f(\mfrak{g},\mfrak{t}_\alpha^-)$ is right exact and the category $\mcal{O}(\mfrak{g})$ has enough projective objects, we may consider the left derived functor
\begin{align*}
LT_\alpha \colon D^b(\mcal{O}(\mfrak{g})) \rarr D^b(\mcal{I}_f(\mfrak{g},\mfrak{t}_\alpha^-))
\end{align*}
of $T_\alpha$ for $\alpha \in \Delta^+$. The following theorem is similar to \cite[Theorem 2.2]{Andersen-Stroppel2003}.
\medskip

\theorem{Let $\alpha \in \Delta^+$. Then we have $L_iT_\alpha = 0$ for all $i >1$. Besides, if $M$ is a $U(\mfrak{s}_\alpha^-)$-free $\mfrak{g}$-module, then we also have $L_1T_\alpha(M) = 0$.}

\proof{Let  $\mcal{A}(\mfrak{g},\mfrak{s}_\alpha^-)$ be the full subcategory of $\mcal{O}(\mfrak{g})$ consisting of $U(\mfrak{s}_\alpha^-)$-free $\mfrak{g}$-modules. Clearly, it contains Verma $\mfrak{g}$-modules and $\mfrak{g}$-modules with a Verma flag, thus all projective objects from $\mcal{O}(\mfrak{g})$. For that reason, the category $\mcal{A}(\mfrak{g},\mfrak{s}_\alpha^-)$ has enough projective objects. Besides, $U(\mfrak{g})_{(f_\alpha)}/U(\mfrak{g}) \simeq \smash{(U(\mfrak{s}_\alpha^-)_{(f_\alpha)}/U(\mfrak{s}_\alpha^-)) \otimes_{U(\mfrak{s}_\alpha^-)} U(\mfrak{g})}$ as left $U(\mfrak{g})$-modules, which gives us that the functor $T_\alpha$ is exact on $\mcal{A}(\mfrak{g},\mfrak{s}_\alpha^-)$. Since every projective object in $\mcal{O}(\mfrak{g})$ is a projective object in $\mcal{A}(\mfrak{g},\mfrak{s}_\alpha^-)$, we get $L_iT_\alpha(M) = 0$ for $i>0$ if $M$ is a $U(\mfrak{s}_\alpha^-)$-free $\mfrak{g}$-module.

Further, let us consider a $\mfrak{g}$-module $M \in \mcal{O}(\mfrak{g})$ and its projective cover $P \rarr M$. Hence, the short exact sequence
$$0 \rarr N \rarr P \rarr M \rarr 0$$ of $\mfrak{g}$-module and the fact that $L_iT_\alpha(P) = 0$ for $i > 0$ gives us $L_iT_\alpha(M) \simeq L_{i-1}T_\alpha(N)$ for $i > 1$. Further, since $U(\mfrak{s}_\alpha^-)$ is a principal ideal domain and $N$ is a $\mfrak{g}$-submodule of a projective object in $\mcal{O}(\mfrak{g})$, i.e.\ $N$ is a $U(\mfrak{s}_\alpha^-)$-free $\mfrak{g}$-module, we have $L_{i-1}T_\alpha(N)=0$ for $i > 1$, which implies $L_iT_\alpha(M) = 0$ for $i>1$.}

\vspace{-2mm}


\subsection{Tensoring with finite-dimensional $\mfrak{g}$-modules}

In this subsection we show that the twisting functors behave well with respect to tensoring with finite-dimensional $\mfrak{g}$-modules, analogously to \cite{Andersen-Stroppel2003}, where this is considered only for simple roots.
\medskip

Let us recall that $U(\mfrak{g})$ is a Hopf algebra with the comultiplication $\Delta \colon U(\mfrak{g}) \rarr U(\mfrak{g}) \otimes_\C U(\mfrak{g})$, the counit $\veps \colon U(\mfrak{g}) \rarr \C$ and the antipode $S \colon U(\mfrak{g}) \rarr U(\mfrak{g})$ given by
\begin{align*}
  \Delta(a) = a \otimes 1 + 1 \otimes a, \qquad \veps(a) = 0, \qquad S(a) = -a
\end{align*}
for $a \in \mfrak{g}$. Since $U(\mfrak{g})_{(f_\alpha)}$ has the structure of a left $\C[f_\alpha^{-1}]$-module, hence also $U(\mfrak{g})_{(f_\alpha)} \otimes_\C U(\mfrak{g})_{(f_\alpha)}$ is a left $\C[f_\alpha^{-1}]$-module and we denote by $U(\mfrak{g})_{(f_\alpha)} \,\smash{\widehat{\otimes}_\C}\, U(\mfrak{g})_{(f_\alpha)}$ its extension to a left $\C[[f_\alpha^{-1}]]$-module, i.e.\ we set
\begin{align*}
  U(\mfrak{g})_{(f_\alpha)} \,\widehat{\otimes}_\C\, U(\mfrak{g})_{(f_\alpha)} = \C[[f_\alpha^{-1}]] \otimes_{\C[f_\alpha^{-1}]}  U(\mfrak{g})_{(f_\alpha)} \otimes_\C U(\mfrak{g})_{(f_\alpha)}.
\end{align*}
There is an obvious extension of the $\C$-algebra structure on $U(\mfrak{g})_{(f_\alpha)} \otimes_\C U(\mfrak{g})_{(f_\alpha)}$ to the completion $U(\mfrak{g})_{(f_\alpha)} \,\smash{\widehat{\otimes}_\C}\, U(\mfrak{g})_{(f_\alpha)}$.
\medskip

Let $\alpha \in \Delta^+$. Then the $\C$-linear mapping
\begin{align*}
  \widetilde{\Delta} \colon U(\mfrak{g})_{(f_\alpha)} &\rarr U(\mfrak{g})_{(f_\alpha)}\, \widehat{\otimes}_\C \, U(\mfrak{g})_{(f_\alpha)}
\end{align*}
given through
\begin{align*}
  \widetilde{\Delta}(f_\alpha^{-n}u) = \bigg(\sum_{k=0}^\infty (-1)^k \binom{n+k-1}{k} f_\alpha^{-n-k} \otimes f_\alpha^k \bigg) \Delta(u),
\end{align*}
where $n \in \N_0$ and $u \in U(\mfrak{g})$, defines a $\C$-algebra homomorphism, see \cite{Andersen-Stroppel2003}. The following theorem is analogous to \cite[Theorem 3.2]{Andersen-Stroppel2003}.
\medskip

\theorem{\label{thm-And}Let $\alpha \in \Delta^+$. Then there is a family $\{\eta_E\}_{E \in \mcal{E}(\mfrak{g})}$ of natural isomorphisms
\begin{align*}
  \eta_E \colon T_\alpha \circ (\,\bullet \otimes_\C E) \rarr (\,\bullet \otimes_\C E) \circ T_\alpha
\end{align*}
of functors such that the following diagrams commute.
\begin{enumerate}[topsep=3pt,itemsep=0pt]
  \item[i)] For $E, F \in \mcal{E}(\mfrak{g})$ we have
  \begin{align*}
    \bfig
    \qtriangle|alr|/->`->`->/<1500,500>[T_\alpha(M \otimes_\C E \otimes_\C F)`T_\alpha(M \otimes_\C E) \otimes_\C F`T_\alpha(M) \otimes_\C E \otimes_\C F.;\eta_F(M \otimes_\C E)`\eta_{E \otimes_\C F}(M)`\eta_E(M) \otimes \id_F]
    \efig
  \end{align*}
  \item[ii)] For $E \in \mcal{E}(\mfrak{g})$ we have
  \begin{align*}
    \bfig
      \square|alrb|/->`->`->`->/<1500,500>[T_\alpha(M \otimes_\C E \otimes_\C E^*)`T_\alpha(M) \otimes_\C E \otimes_\C E^*`T_\alpha(M \otimes_\C \C)`T_\alpha(M) \otimes_\C \C,;\eta_{E \otimes_\C E^*}(M)`T_\alpha(\id_M \otimes i_E)`\id_{T_\alpha(M)} \otimes i_E`\eta_\C(M)]
    \efig
  \end{align*}
  where $i_E \colon \C \rarr E \otimes_\C E^*$ is given by $1 \mapsto \sum_{i=1}^d e_i \otimes e_i^*$ for a fixed basis $\{e_i\}_{1 \leq i \leq d}$ of $E$ with the dual basis $\{e_i^*\}_{1 \leq i \leq d}$ of $E^*$.
\end{enumerate}}

\remark{Theorem \ref{thm-And} implies that for any $\alpha \in \Delta^+$ the twisting functor $T_\alpha$ commutes with the translation functors.}

\vspace{-2mm}


\subsection{Twisting of generalized Verma modules}

In Section \ref{subsec:Twisting functor} we introduced the twisting functor $T_\alpha$ and the partial Zuckerman functor $S_\alpha$ for $\alpha \in \Delta^+$ as endofunctors of the category $\mcal{M}(\mfrak{g})$. Now, to construct new simple weight $\mfrak{g}$-modules with infinite-dimensional weight spaces and finite $\Gamma_\alpha$-multiplicities we
 restrict both functors to the  subcategory $\mcal{O}(\mfrak{g})$. Therefore, we have
\begin{align*}
  T_\alpha \colon \mcal{O}(\mfrak{g}) \rarr \mcal{I}_f(\mfrak{g},\mfrak{t}_\alpha^-) \qquad \text{and} \qquad  S_\alpha \colon \mcal{O}(\mfrak{g}) \rarr \mcal{I}_f(\mfrak{g},\mfrak{t}_\alpha^-),
\end{align*}
which follows by Theorem \ref{thm:finitely generated}.

In addition to the standard Borel subalgebra of $\mfrak{g}$ we  also consider the standard parabolic subalgebras of $\mfrak{g}$. For a subset $\Sigma$ of $\Pi$  denote by $\Delta_\Sigma$ the root subsystem in $\mfrak{h}^*$ generated by $\Sigma$. Then the standard parabolic subalgebra $\mfrak{p}_\Sigma$ of $\mfrak{g}$ associated to $\Sigma$ is defined as $\mfrak{p}_\Sigma = \mfrak{l}_\Sigma \oplus \mfrak{u}_\Sigma$ with the nilradical $\mfrak{u}_\Sigma$ and the opposite nilradical $\widebar{\mfrak{u}}_\Sigma$ given by
\begin{align*}
  \mfrak{u}_\Sigma = \bigoplus_{\alpha \in \Delta^+ \setminus \Delta_\Sigma} \mfrak{g}_\alpha \qquad \text{and} \qquad \widebar{\mfrak{u}}_\Sigma = \bigoplus_{\alpha \in \Delta^+ \setminus \Delta_\Sigma} \mfrak{g}_{-\alpha}
\end{align*}
and with the Levi subalgebra $\mfrak{l}_\Sigma$ defined by
\begin{align*}
  \mfrak{l}_\Sigma = \mfrak{h} \oplus \bigoplus_{\alpha \in \Delta_\Sigma} \mfrak{g}_\alpha.
\end{align*}
Moreover, we have  the corresponding triangular decomposition
\begin{align*}
  \mfrak{g}= \widebar{\mfrak{u}}_\Sigma \oplus \mfrak{l}_\Sigma \oplus \mfrak{u}_\Sigma
 \end{align*}
of the Lie algebra $\mfrak{g}$. Note that if $\Sigma =\emptyset$ then $\mfrak{p}_\Sigma = \mfrak{b}$ and if $\Sigma = \Pi$ then $\mfrak{p}_\Sigma = \mfrak{g}$.
\medskip

Let $\mfrak{p}=\mfrak{l} \oplus \mfrak{u}$ be the standard parabolic subalgebra of $\mfrak{g}$ associated to a subset $\Sigma$ of $\Pi$. Denote
\begin{align*}
  \Delta^+_\mfrak{u} = \{\alpha \in \Delta^+;\, \mfrak{g}_\alpha \subset \mfrak{u}\}, \qquad \Delta^+_\mfrak{l} = \{\alpha \in \Delta^+;\, \mfrak{g}_\alpha \subset \mfrak{l}\}
\end{align*}
and set
\begin{align*}
  \Lambda^+(\mfrak{p}) = \{\lambda \in \mfrak{h}^*;\, (\forall \alpha \in \Sigma)\, \lambda(h_\alpha) \in \N_0\}.
\end{align*}
The elements of $\Lambda^+(\mfrak{p})$ are called $\mfrak{p}$-dominant and $\mfrak{p}$-algebraically integral weights.

Besides, we denote by $\sigma_\lambda \colon \mfrak{p} \rarr \End \mathbb{F}_\lambda$ the simple finite-dimensional $\mfrak{p}$-module with highest weight $\lambda \in \Lambda^+(\mfrak{p})$. Let us note that the nilradical $\mfrak{u}$ of $\mfrak{p}$ acts trivially on $\mathbb{F}_\lambda$.
\medskip

\definition{Let $\lambda \in \Lambda^+(\mfrak{p})$. The generalized Verma $\mfrak{g}$-module $M^\mfrak{g}_\mfrak{p}(\lambda)$ is the induced module
\begin{align}
  M^\mfrak{g}_\mfrak{p}(\lambda) = \Ind^\mfrak{g}_\mfrak{p} \mathbb{F}_\lambda = U(\mfrak{g}) \otimes_{U(\mfrak{p})} \mathbb{F}_\lambda \simeq U(\widebar{\mfrak{u}}) \otimes_\C \mathbb{F}_\lambda,
\end{align}
where the last isomorphism of $U(\widebar{\mfrak{u}})$-modules follows from the Poincaré--Birkhoff--Witt theorem.}

Now, we are in the position to define our main objects -- $\mfrak{g}$-modules which are twistings of generalized Verma modules. Later we will give their explicit realization.
\medskip

\definition{Let $\lambda \in \Lambda^+(\mfrak{p})$ and $\alpha \in \Delta^+_\mfrak{u}$. Then the $\alpha$-Gelfand--Tsetlin $\mfrak{g}$-module $W^\mfrak{g}_\mfrak{p}(\lambda,\alpha)$ is defined by
\begin{align}
  W^\mfrak{g}_\mfrak{p}(\lambda,\alpha) = T_\alpha(M^\mfrak{g}_\mfrak{p}(\lambda)),
\end{align}
where $T_\alpha \colon \mcal{O}(\mfrak{g}) \rarr \mcal{I}_f(\mfrak{g},\mfrak{t}_\alpha^-)$ is the twisting functor.}

Let us note that for $\lambda \in \Lambda^+(\mfrak{p})$ and $\alpha \in \Delta^+_\mfrak{l}$ we have  $T_\alpha(M^\mfrak{g}_\mfrak{p}(\lambda))=0$ by Proposition \ref{prop:locally nilpotent}. This explains the restriction $\alpha \in \Delta^+_\mfrak{u}$ in the definition above.
\medskip

By using the triangular decomposition $\mfrak{g} = \widebar{\mfrak{u}} \oplus \mfrak{l} \oplus \mfrak{u}$ and the Poincaré--Birkhoff--Witt theorem we get an isomorphism
\begin{align*}
  U(\mfrak{g})_{(f_\alpha)} \simeq U(\widebar{\mfrak{u}})_{(f_\alpha)} \otimes_\C U(\mfrak{p})
\end{align*}
of left $U(\widebar{\mfrak{u}})$-modules for $\alpha \in \Delta^+$. Hence, we may write
\begin{align*}
  U(\mfrak{g})_{(f_\alpha)} \otimes_{U(\mfrak{g})} M^\mfrak{g}_\mfrak{p}(\lambda) \simeq U(\mfrak{g})_{(f_\alpha)} \otimes_{U(\mfrak{g})} U(\mfrak{g}) \otimes_{U(\mfrak{p})} \mathbb{F}_\lambda \simeq U(\widebar{\mfrak{u}})_{(f_\alpha)} \otimes_\C \mathbb{F}_\lambda
\end{align*}
for $\lambda \in \Lambda^+(\mfrak{p})$ and $\alpha \in \Delta^+_\mfrak{u}$, which gives us an isomorphism
\begin{align*}
  W^\mfrak{g}_\mfrak{p}(\lambda,\alpha) \simeq (U(\widebar{\mfrak{u}})_{(f_\alpha)}/U(\widebar{\mfrak{u}})) \otimes_\C \mathbb{F}_\lambda
\end{align*}
of $U(\widebar{\mfrak{u}})$-modules. Besides, we have that
\begin{align}
  M^\mfrak{g}_\mfrak{p}(\lambda) \simeq U(\widebar{\mfrak{u}}) \otimes_\C \mathbb{F}_\lambda \qquad \text{and} \qquad W^\mfrak{g}_\mfrak{p}(\lambda,\alpha) \simeq (U(\widebar{\mfrak{u}})_{(f_\alpha)}/U(\widebar{\mfrak{u}})) \otimes_\C \mathbb{F}_\lambda
\end{align}
are also isomorphisms of $U(\mfrak{l})$-modules for the adjoint action of $U(\mfrak{l})$ on $U(\widebar{\mfrak{u}})$ and $U(\widebar{\mfrak{u}})_{(f_\alpha)}$.

For later use, we also want to clarify the relation between $W^\mfrak{g}_\mfrak{b}(\lambda,\alpha)$ and $W^\mfrak{g}_\mfrak{p}(\lambda,\alpha)$ for $\lambda \in \Lambda^+(\mfrak{p})$ and $\alpha \in \Delta^+_\mfrak{u}$. Since we have a canonical surjective homomorphism
\begin{align*}
  M^\mfrak{g}_\mfrak{b}(\lambda) \rarr M^\mfrak{g}_\mfrak{p}(\lambda)
\end{align*}
of generalized Verma modules, by applying the twisting functor $T_\alpha$, which is right exact, we get a surjective homomorphism
\begin{align*}
  W^\mfrak{g}_\mfrak{b}(\lambda,\alpha) \rarr W^\mfrak{g}_\mfrak{p}(\lambda,\alpha)
\end{align*}
of $\alpha$-Gelfand--Tsetlin modules. Therefore, we obtain that $W^\mfrak{g}_\mfrak{p}(\lambda,\alpha)$ is a quotient of $W^\mfrak{g}_\mfrak{b}(\lambda,\alpha)$ for $\lambda \in \Lambda^+(\mfrak{p})$ and $\alpha \in \Delta^+_\mfrak{u}$.


\section{$\alpha$-Gelfand--Tsetlin modules $W^\mfrak{g}_\mfrak{p}(\lambda,\alpha)$}
\label{sec:Realizations}

In this section we discuss the properties of $\alpha$-Gelfand--Tsetlin module $W^\mfrak{g}_\mfrak{p}(\lambda,\alpha)$ for $\lambda \in \Lambda^+(\mfrak{p})$ and $\alpha \in \Delta^+_\mfrak{u}$ and show that they belong to the category $\mcal{H}_{\rm fin}(\mfrak{g}, \Gamma_\alpha)$. We also give their explicit realization.


\subsection{Basic properties}

We denote by $\Z^{\Delta^+}$\! and $\smash{\Z^{\Delta^+_\mfrak{u}}}$ the set of all functions from $\Delta^+$ to $\Z$ and from $\smash{\Delta^+_\mfrak{u}}$ to $\Z$, respectively. Since  $\Delta^+_\mfrak{u} \subset \Delta^+$, an element of $\Z^{\Delta^+_\mfrak{u}}$ will be also regarded as an element of $\Z^{\Delta^+}$\! extended by $0$ on $\Delta^+ \setminus \Delta^+_\mfrak{u}$. A similar notation is introduced for $\smash{\N_0^{\Delta^+}}$\! and $\smash{{\N_0}^{\!\!\!\Delta^+_\mfrak{u}}}$.
\medskip

Since any positive root $\gamma \in \Delta^+$ can be expressed as
\begin{align}
  \gamma = \sum_{i=1}^r m_{\gamma,\alpha_i} \alpha_i, \label{eq:simple root decomposition}
\end{align}
where $m_{\gamma,\alpha_i} \in \N_0$ for $i=1,2,\dots,r$ and $\Pi=\{\alpha_1, \alpha_2, \dots, \alpha_r\}$, we define $t_\gamma^\alpha \in \Z^{\Delta^+}$\! for $\gamma \in \Delta^+ \setminus \Pi$ and $\alpha \in \Delta^+$ by
\begin{align*}
  t_\gamma^\alpha(\beta) = \begin{cases}
    -(-1)^{\delta_{\alpha, \alpha_i}}m_{\gamma,\alpha_i} &  \text{for $\beta=\alpha_i$}, \\
    (-1)^{\delta_{\alpha,\gamma}} & \text{for $\beta=\gamma$}, \\
    0 & \text{for $\beta \neq \gamma$ and $\beta \notin \Pi$},
  \end{cases}
\end{align*}
and the subset $\Lambda_+^\alpha$ of $\Z^{\Delta^+}$\! by
\begin{align}
  \Lambda_+^\alpha = \{{\textstyle \sum_{\gamma \in \Delta^+ \setminus \Pi} n_\gamma t_\gamma^\alpha};\, n_\gamma \in \N_0\ \text{for all $\gamma \in \Delta^+ \setminus \Pi$}\}.
\end{align}
Besides, we introduce a weight $\mu_{a,\alpha} \in \mfrak{h}^*$ by
\begin{align}
  \mu_{a,\alpha} = -\sum_{\gamma \in \Delta^+} (-1)^{\delta_{\alpha,\gamma}} a_\gamma \gamma + \alpha \label{eq:weights}
\end{align}
for $\alpha \in \Delta^+$ and $a \in \Z^{\Delta^+}$\!.
\medskip

\lemma{\label{lem:mu weights}
Let $\alpha \in \Delta^+_\mfrak{u}$ and $a,b \in \Z^{\Delta^+_\mfrak{u}}$. Then $\mu_{a,\alpha} = \mu_{b,\alpha}$ if and only if
\begin{align}
  b = a + \sum_{\gamma \in \Delta^+_\mfrak{u} \setminus \Pi} n_\gamma t_\gamma^\alpha,
\end{align}
where $n_\gamma \in \Z$ for all $\gamma \in \Delta^+_\mfrak{u} \setminus \Pi$.}

\proof{Let $\alpha \in \Delta^+_\mfrak{u}$ and $a,b \in \Z^{\Delta^+_\mfrak{u}}$. If
\begin{align*}
b=a+ \sum_{\gamma \in \Delta^+_\mfrak{u} \setminus \Pi} n_\gamma t_\gamma^\alpha,
\end{align*}
where $n_\gamma \in \Z$ for $\gamma \in \Delta^+_\mfrak{u} \setminus \Pi$, then we easily get $\mu_{a,\alpha}=\mu_{b,\alpha}$. On the other hand, let us assume that $\mu_{a,\alpha}=\mu_{b,\alpha}$. Then we set $n_\gamma = (-1)^{\delta_{\alpha,\gamma}}(b_\gamma-a_\gamma)$ for $\gamma \in \Delta^+_\mfrak{u} \setminus \Pi$ and define $c=a + \sum_{\gamma \in \Delta^+_\mfrak{u} \setminus \Pi} n_\gamma t_\gamma^\alpha$. Hence, we have $c_\gamma = a_\gamma + (-1)^{\delta_{\alpha,\gamma}}n_\gamma = b_\gamma$ for $\gamma \in \Delta^+_\mfrak{u} \setminus \Pi$ and
\begin{align*}
c_{\alpha_i} = a_{\alpha_i} - \sum_{\gamma \in \Delta^+_\mfrak{u} \setminus \Pi} (-1)^{\delta_{\alpha,\alpha_i}}n_\gamma m_{\gamma,\alpha_i}
\end{align*}
for $i=1,2,\dots,r$. Further, since $\mu_{a,\alpha}=\mu_{b,\alpha}$, we may write
\begin{align*}
  \mu_{b,\alpha}-\mu_{a,\alpha} &= \sum_{\gamma \in \Delta^+_\mfrak{u}} (-1)^{\delta_{\alpha,\gamma}}(a_\gamma-b_\gamma) \gamma  \\
  &= \sum_{i=1}^r (-1)^{\delta_{\alpha,\alpha_i}} (a_{\alpha_i}-b_{\alpha_i})\alpha_i + \sum_{i=1}^r \sum_{\gamma \in \Delta^+_\mfrak{u} \setminus \Pi} (-1)^{\delta_{\alpha,\gamma}}(a_\gamma - b_\gamma)m_{\gamma,\alpha_i}\alpha_i \\
  &= \sum_{i=1}^r (-1)^{\delta_{\alpha,\alpha_i}}\!\Big(a_{\alpha_i}-b_{\alpha_i} - \sum_{\gamma \in \Delta^+_\mfrak{u} \setminus \Pi} (-1)^{\delta_{\alpha,\alpha_i}} n_\gamma m_{\gamma,\alpha_i} \Big)\alpha_i =0,
\end{align*}
where we used \eqref{eq:simple root decomposition}. As the set $\{\alpha_1,\alpha_2,\dots,\alpha_r\}$ forms a basis of $\mfrak{h}^*$, we get
\begin{align*}
 a_{\alpha_i}-b_{\alpha_i} - \sum_{\gamma \in \Delta^+_\mfrak{u} \setminus \Pi} (-1)^{\delta_{\alpha,\alpha_i}} n_\gamma m_{\gamma,\alpha_i}=0,
\end{align*}
which implies that $c_{\alpha_i}=b_{\alpha_i}$ for $i=1,2,\dots,r$. Hence, we have $c_\gamma= b_\gamma$ for all $\gamma \in \Delta^+_\mfrak{u}$ and we are done.}

Let $\{f_\alpha,f_{\gamma_1},\dots,f_{\gamma_n}\}$ be a root basis of the opposite nilradical $\widebar{\mfrak{u}}$, where $n = \dim \widebar{\mfrak{u}}-1$, $\gamma_i \in \Delta^+_\mfrak{u}$ and $f_{\gamma_i} \in \mfrak{g}_{-\gamma_i}$ for $i=1,2,\dots,n$. Then as a consequence of the Poincaré--Birkhoff--Witt theorem we obtain that the subset $\{u_{a,\alpha};\, a \in \smash{{\N_0}^{\!\!\!\Delta^+_\mfrak{u}}}\}$ of $U(\widebar{\mfrak{u}})_{(f_\alpha)}/U(\widebar{\mfrak{u}})$, where
\begin{align*}
  u_{a,\alpha} = f_\alpha^{-a_\alpha-1} f_{\gamma_1}^{a_{\gamma_1}} \dotsm f_{\gamma_n}^{a_{\gamma_n}}
\end{align*}
for $a \in \smash{{\N_0}^{\!\!\!\Delta^+_\mfrak{u}}}$, forms a basis of $U(\widebar{\mfrak{u}})_{(f_\alpha)}/U(\widebar{\mfrak{u}})$ for $\alpha \in \Delta^+_\mfrak{u}$. Moreover, if $v \in \mathbb{F}_\lambda$ is a weight vector with weight $\mu_v \in \mfrak{h}^*$, then $u_{a,\alpha} \otimes v \in W^\mfrak{g}_\mfrak{p}(\lambda,\alpha)$ is a weight vector with weight $\mu_v + \mu_{a,\alpha}$  by Lemma \ref{lem:sl2 action}.
\medskip

\proposition{\label{prop:weight spaces}
Let $\lambda  \in \Lambda^+(\mfrak{p})$ and $\alpha \in \Delta^+_\mfrak{u}$. Then all weight spaces of $W^\mfrak{g}_\mfrak{p}(\lambda,\alpha)$ are finite dimensional if $\alpha \in \Delta^+_\mfrak{u} \cap \Pi$ and infinite dimensional if $\alpha \in \Delta^+_\mfrak{u} \setminus \Pi$.}

\proof{Since $W^\mfrak{g}_\mfrak{p}(\lambda,\alpha)$  is isomorphic to $(U(\widebar{\mfrak{u}})_{(f_\alpha)}/U(\widebar{\mfrak{u}})) \otimes_\C \mathbb{F}_\lambda$ as an $\mfrak{l}$-module  for $\lambda \in \Lambda^+(\mfrak{p})$ and $\alpha \in \Delta^+_\mfrak{u}$, and the $\mfrak{g}$-module $\mathbb{F}_\lambda$ is finite dimensional, it is sufficient to show that all weight spaces of $U(\widebar{\mfrak{u}})_{(f_\alpha)}/U(\widebar{\mfrak{u}})$ are finite dimensional if $\alpha \in \Delta^+_\mfrak{u} \cap \Pi$ and infinite dimensional if $\alpha \in \Delta^+_\mfrak{u} \setminus \Pi$.

From the discussion above it follows that all weights of $U(\widebar{\mfrak{u}})_{(f_\alpha)}/U(\widebar{\mfrak{u}})$ are of the form $\mu_{a,\alpha}$ for some $a \in \Z^\Pi$. If we assume that $u_{b,\alpha}$ for $b \in \smash{{\N_0}^{\!\!\!\Delta^+_\mfrak{u}}}$\! is a weight vector with weight $\mu_{a,\alpha}$, then we obtain $\mu_{a,\alpha}=\mu_{b,\alpha}$, which together with Lemma \ref{lem:mu weights} implies
\begin{align*}
  b = a + \sum_{\gamma\in \Delta^+ \setminus \Pi} n_\gamma t_\gamma^\alpha,
\end{align*}
where $n_\gamma \in \Z$ for all $\gamma \in \Delta^+ \setminus \Pi$. As $\{u_{b,\alpha};\, b \in \smash{{\N_0}^{\!\!\!\Delta^+_\mfrak{u}}}\}$ is a basis of $U(\widebar{\mfrak{u}})_{(f_\alpha)}/U(\widebar{\mfrak{u}})$, we immediately get
\begin{align*}
   \dim\, (U(\widebar{\mfrak{u}})_{(f_\alpha)}/U(\widebar{\mfrak{u}}))_{\mu_{a,\alpha}} = \natural \{t \in \Lambda_+^\alpha;\, a + t \in \smash{{\N_0}^{\!\!\!\Delta^+_\mfrak{u}}}\}.
\end{align*}
Let us assume that $t = \smash{\sum_{\gamma \in \Delta^+ \setminus \Pi}}\, n_\gamma t^\alpha_\gamma \in \Lambda_+^\alpha$ satisfies $a + t \in \smash{{\N_0}^{\!\!\!\Delta^+_\mfrak{u}}}$. Then for $b=a+t$ we have $b_\gamma = a_\gamma +(-1)^{\delta_{\alpha,\gamma}}n_\gamma = (-1)^{\delta_{\alpha,\gamma}}n_\gamma$ for $\gamma \in \Delta^+ \setminus \Pi$ and $b_{\alpha_i}= a_{\alpha_i} - \smash{\sum_{\gamma \in \Delta^+ \setminus \Pi}}\, (-1)^{\delta_{\alpha,\alpha_i}} n_\gamma m_{\gamma,\alpha_i}$ for $i=1,2,\dots,r$.

If $\alpha \in \Delta^+_\mfrak{u} \cap \Pi$, then the condition $b_\gamma \in \N_0$ for $\gamma \in \Delta^+$ implies $n_\gamma \in \N_0$ for $\gamma \in \Delta^+ \setminus \Pi$. Further, for $i=1, 2, \dots, r$ satisfying $\alpha_i \neq \alpha$ we may write
\begin{align*}
  0 \leq b_{\alpha_i} = a_{\alpha_i} - \sum_{\gamma \in \Delta^+ \setminus \Pi} n_\gamma m_{\gamma,\alpha_i} \leq a_{\alpha_i} - n_\gamma m_{\gamma,\alpha_i}.
\end{align*}
Moreover, for each $\gamma \in \Delta^+ \setminus \Pi$ there exists $i \in \{1,2,\dots,r\}$ such that $\alpha_i \neq \alpha$ and $m_{\gamma,\alpha_i} \neq 0$, which implies that $\dim\, (U(\widebar{\mfrak{u}})_{(f_\alpha)}/U(\widebar{\mfrak{u}}))_{\mu_{a,\alpha}} < \infty$.

On the other hand, if $\alpha \in \Delta^+_\mfrak{u} \setminus \Pi$, then $a + t - nt_\alpha^\alpha \in \smash{{\N_0}^{\!\!\!\Delta^+_\mfrak{u}}}$ for all $n \in \N_0$ provided $a + t \in \smash{{\N_0}^{\!\!\!\Delta^+_\mfrak{u}}}$, which implies that $\dim\, (U(\widebar{\mfrak{u}})_{(f_\alpha)}/U(\widebar{\mfrak{u}}))_{\mu_{a,\alpha}} = \infty$. This finishes the proof.}

The next theorem generalizes \cite[Theorem 2.9, Theorem 2.11]{Futorny-Krizka2017} for any parabolic subalgebra $\mfrak{p}$ of $\mfrak{g}$ and for a positive root $\alpha \in \Delta^+_\mfrak{u}$.
\medskip

\theorem{\label{thm:GT modules}
Let $\lambda \in \Lambda^+(\mfrak{p})$ and $\alpha \in \Delta^+_\mfrak{u}$. Then $W^\mfrak{g}_\mfrak{p}(\lambda,\alpha)$ is a locally $\mfrak{t}_\alpha^-$-finite cyclic weight $\mfrak{g}$-module with finite $\Gamma_\alpha$-multiplicities, that is
$W^\mfrak{g}_\mfrak{p}(\lambda,\alpha)\in\mcal{H}_{\rm fin}(\mfrak{g},\Gamma_\alpha)\cap \mcal{I}_f(\mfrak{g}, \mfrak{t}_\alpha^-)$.}

\proof{It is an immediate consequence of Theorem \ref{thm:highest weight modules condition} (ii), since we have $W^\mfrak{g}_\mfrak{p}(\lambda,\alpha) = T_\alpha(M^\mfrak{g}_\mfrak{p}(\lambda))$ for $\lambda \in \Lambda^+(\mfrak{p})$ and $\alpha \in \Delta^+_\mfrak{u}$.}

\remark{It was shown in \cite[Theorem 2.16]{Futorny-Krizka2017} that the $\mfrak{g}$-module $W^\mfrak{g}_\mfrak{b}(\lambda,\theta)$ is simple generically for $\mfrak{g}=\mfrak{sl}(3)$, where $\theta$ is the maximal root of $\mfrak{g}$. Moreover, in this case $\Gamma_\theta$ is diagonalizable on $W^\mfrak{g}_\mfrak{b}(\lambda,\theta)$. Therefore, it is natural to expect a similar behaviour from $W^\mfrak{g}_\mfrak{p}(\lambda,\alpha)$ for $\lambda \in \Lambda^+(\mfrak{p})$ and $\alpha \in \Delta^+_\mfrak{u}$ in general. }

The $\mfrak{g}$-module $W^\mfrak{g}_\mfrak{p}(\lambda,\alpha)$ for $\lambda \in \Lambda^+(\mfrak{p})$ and $\alpha \in \Delta^+_\mfrak{u}$  is an $\alpha$-Gelfand--Tsetlin module in
the category $\mcal{I}_{\rm fin}(\mfrak{g},\mfrak{s}_\alpha^-)$.
 On the other hand, by means of the natural duality we can analogously construct certain $\alpha$-Gelfand--Tsetlin modules in the category $\mcal{I}_{\rm fin}(\mfrak{g}, \mfrak{s}_\alpha^+)$.
\medskip

Furthermore, since the action of the center $\mfrak{z}_\mfrak{g}$ of $U(\mfrak{g})$ on the generalized Verma module $M^\mfrak{g}_\mfrak{p}(\lambda)$ for $\lambda \in \Lambda^+(\mfrak{p})$ is then given by
\begin{align}
  zv=\chi_{\lambda+\rho}(z)v \label{eq:center action Verma}
\end{align}
for all $z \in \mfrak{z}_\mfrak{g}$ and $v \in M^\mfrak{g}_\mfrak{p}(\lambda)$,
in other words $M^\mfrak{g}_\mfrak{p}(\lambda)$ is a $\mfrak{g}$-module with central character $\chi_{\lambda+\rho}$, it follows immediately from definition that $W^\mfrak{g}_\mfrak{p}(\lambda,\alpha)$ for $\lambda \in \Lambda^+(\mfrak{p})$ and $\alpha \in \Delta^+_\mfrak{u}$ is a $\mfrak{g}$-module with the same central character $\chi_{\lambda+\rho}$.
\medskip

A comparison of basic characteristics of generalized Verma modules $M^\mfrak{g}_\mfrak{p}(\lambda)$ and $\alpha$-Gelfand--Tsetlin modules $W^\mfrak{g}_\mfrak{p}(\lambda,\alpha)$ for $\lambda \in \Lambda^+(\mfrak{p})$ and $\alpha \in \Delta^+_\mfrak{u}$ is given in Table \ref{tab:comparision}.


\subsection{Geometric realization of $W^\mfrak{g}_\mfrak{b}(\lambda,\alpha)$}

We describe the geometric realization of $\alpha$-Gelfand--Tsetlin modules $W^\mfrak{g}_\mfrak{b}(\lambda,\alpha)$ for $\lambda \in \mfrak{h}^*$ and $\alpha \in \Delta^+$
using the theory of algebraic $\mcal{D}$-modules on flag varieties without going into details. We are going to address the geometric realization of more general families of Gelfand--Tsetlin modules and the corresponding geometric induction in the subsequent paper.
\medskip

Let $G$ be a connected semisimple algebraic group over $\C$, $H$ be a maximal torus of $G$ and $B$ be a Borel subgroup of $G$ containing $H$ with the unipotent radical $N$ and the opposite unipotent radical $\widebar{N}$. We denote by $\mfrak{g}$, $\mfrak{n}$, $\widebar{\mfrak{n}}$ and $\mfrak{h}$ the Lie algebras of $G$, $N$, $\widebar{N}$ and $H$, respectively. Then
\begin{align*}
  X = G/B
\end{align*}
is a smooth algebraic variety, the \emph{flag variety} of $G$. Besides, we have the canonical $G$-equivariant projection
\begin{align*}
  p \colon G \rarr G/B.
\end{align*}
Following \cite{Kashiwara1989}, for any $\lambda \in \mfrak{h}^*$ there exists a $G$-equivariant sheaf of rings of twisted differential operators $\mcal{D}_X^\lambda$ on $X$. Let us note that $\mcal{D}_X^{-\rho}$ is the usual sheaf of rings of differential operators on $X$.
Since $\mcal{D}_X^\lambda$ is $G$-equivariant, we have a Lie algebra homomorphism $\mfrak{g} \rarr \Gamma(X,\mcal{D}_X^\lambda)$, which extends to a homomorphism
\begin{align}
  \Phi_\lambda \colon U(\mfrak{g}) \rarr \Gamma(X,\mcal{D}_X^\lambda) \label{eq:Phi homomorphism}
\end{align}
of $\C$-algebras. Hence, for any $\mcal{D}_X^\lambda$-module $\mcal{M}$ the vector space $\Gamma(X,\mcal{M})$ of global sections of $\mcal{M}$ has  a natural $\mfrak{g}$-module structure.
\medskip

Let us denote by $\Mod_c(\mcal{D}_X^\lambda)$ the category of coherent $\mcal{D}_X^\lambda$-modules and by $\mcal{M}_f(\mfrak{g},\chi_\lambda)$ the full subcategory of $\mcal{M}(\mfrak{g})$ consisting of finitely generated $\mfrak{g}$-modules on which the center $\mfrak{z}_\mfrak{g}$ of $U(\mfrak{g})$ acts via the central character $\chi_\lambda$. The following remarkable theorem is due to Beilinson and Bernstein, see \cite{Beilinson-Bernstein1981}, \cite{Kashiwara1989}.
\medskip

\theorem{\label{thm:BB correspondence}
Let $\lambda \in \mfrak{h}^*$ be anti-dominant and regular, i.e.\ $\lambda(h_\alpha) \notin \N_0$ for $\alpha \in \Delta^+$. Then the functor of global sections
\begin{align*}
  \Gamma(X,\bullet\,) \colon  \Mod_c(\mcal{D}_X^\lambda) \rarr \mcal{M}_f(\mfrak{g},\chi_\lambda)
\end{align*}
induces an equivalence of abelian categories. The inverse is given through the localization functor $\Delta(M) = \mcal{D}_X^\lambda \otimes_{U(\mfrak{g})}M$ for $M \in \mcal{M}_f(\mfrak{g},\chi_\lambda)$.}

Since the $\alpha$-Gelfand--Tsetlin module $W^\mfrak{g}_\mfrak{b}(\lambda-\rho,\alpha)$ for $\lambda \in \mfrak{h}^*$ and $\alpha \in \Delta^+$ is an object of the category $\mcal{M}_f(\mfrak{g},\chi_\lambda)$, it corresponds to the coherent $\mcal{D}_X^\lambda$-module on the flag variety $X$ by Theorem \ref{thm:BB correspondence}. Our goal is to describe this coherent $\mcal{D}_X^\lambda$-module explicitly.
\medskip

First of all we recall the geometric realization of Verma modules $M^\mfrak{g}_\mfrak{b}(\lambda-\rho)$ for $\lambda \in \mfrak{h}^*$. Let $X_w$ for $w \in W \simeq N_G(H)/H$ be the $N$-orbit in $X$ defined by
\begin{align*}
  X_w = N\dot{w}B/B,
\end{align*}
where $\dot{w} \in N_G(H)$ is a representative of $w$. Let us note that $X_w$ is called a \emph{Schubert cell}. Further, let us denote by
\begin{align*}
  i_w \colon X_w \rarr X
\end{align*}
the embedding of $X_w$ into $X$ for $w \in W$. Since $i_w$ is an $N$-equivariant mapping, the pull-back $i_w^\sharp \mcal{D}_X^\lambda$ for $\lambda \in \mfrak{h}^*$ is an $N$-equivariant sheaf of rings of twisted differential operators on $X_w$, which is isomorphic to the sheaf of rings of differential operators $\mcal{D}_{X_w}$\! on $X_w$. As the structure sheaf $\mcal{O}_{X_w}$\! of $X_w$ is a $\mcal{D}_{X_w}$-module, we may consider its direct image
\begin{align}
  \mcal{L}_\lambda(X_w,\mcal{O}_{X_w}) = i_{w*}(\mcal{D}_{X \larr X_w}^\lambda \!\otimes_{i_w^\sharp \mcal{D}_X^\lambda} \mcal{O}_{X_w}),
\end{align}
where the $(i_w^{-1} \mcal{D}_X^\lambda,i_w^\sharp \mcal{D}_X^\lambda)$-bimodule $\mcal{D}_{X \larr X_w}^\lambda$\! is the so-called transfer bimodule for $i_w$. Then for the unit element $e \in W$ we have
\begin{align}
  \Gamma(X,\mcal{L}_\lambda(X_e,\mcal{O}_{X_e})) \simeq M^\mfrak{g}_\mfrak{b}(\lambda-\rho)
\end{align}
as $\mfrak{g}$-modules for $\lambda \in \mfrak{h}^*$, see e.g.\ \cite{Krizka-Somberg2017}.
\medskip

To obtain the geometric realization of $\alpha$-Gelfand--Tsetlin modules $W^\mfrak{g}_\mfrak{b}(\lambda-\rho,\alpha)$ for $\lambda \in \mfrak{h}^*$ and $\alpha \in \Delta^+$, we take into account different embeddings into $X$. Since $\mfrak{s}_\alpha$ is a semisimple Lie subalgebra of $\mfrak{g}$ for any $\alpha \in \Delta^+$, there is a unique closed connected semisimple algebraic subgroup $G_\alpha$ of $G$ with the Lie algebra $\mfrak{s}_\alpha$. Moreover, we have that $H_\alpha = H \cap G_\alpha$ is a maximal torus of $G_\alpha$ and $B_\alpha = B \cap G_\alpha$ is a Borel subgroup of $G_\alpha$ with the unipotent radical $N_\alpha = N \cap G_\alpha$ and the opposite unipotent radical $\widebar{N}_\alpha = \widebar{N} \cap G_\alpha$. Let us note that $\mfrak{s}_\alpha^+$ and $\mfrak{s}_\alpha^-$ are Lie algebras of $N_\alpha$ and $\widebar{N}_\alpha$, respectively. Let us consider the closed embedding
\begin{align*}
  i_\alpha \colon X_\alpha \rarr X
\end{align*}
of the flag variety $X_\alpha = G_\alpha/B_\alpha$ into $X$. Moreover, since $i_\alpha$ is a $G_\alpha$-equivariant mapping, the pull-back $i_\alpha^\sharp \mcal{D}_X^\lambda$ is a $G_\alpha$-equivariant sheaf of rings of twisted differential operators on $X_\alpha$, which is isomorphic to \smash{$\mcal{D}_{X_\alpha}^{\lambda_\alpha}$} with $\lambda_\alpha= \big(\lambda + \rho - {1 \over 2}\alpha\big)|_{\mfrak{h} \cap \mfrak{s}_\alpha}$. Hence, for $\alpha \in \Delta^+$ we may introduce the geometric induction functor
\begin{align*}
  \GInd_\alpha^\lambda \colon \mcal{M}_f(\mfrak{s}_\alpha,\chi_{\lambda_\alpha}) \rarr \mcal{M}_f(\mfrak{g},\chi_\lambda)
\end{align*}
provided $\lambda \in \mfrak{h}^*$ and $\lambda_\alpha \in (\mfrak{h} \cap \mfrak{s}_\alpha)^*$ are anti-dominant and regular by
\begin{align}
  \GInd_\alpha^\lambda(M) = \Gamma(X,i_{\alpha*}(\mcal{D}^\lambda_{X \larr X_\alpha} \! \otimes_{i_\alpha^\sharp \mcal{D}_X^\lambda} \Delta(M))),
\end{align}
where the $(i_\alpha^{-1} \mcal{D}_X^\lambda,i_\alpha^\sharp \mcal{D}_X^\lambda)$-bimodule $\mcal{D}_{X \larr X_\alpha}^\lambda$ is the transfer bimodule for $i_\alpha$.
\medskip

Let us consider an open subset $U_e=p(\widebar{N})$ of $X$ called the big cell. Then $\mcal{D}_X^\lambda|_{U_e}$ and $\mcal{D}_X|_{U_e}$ are isomorphic as sheaves of rings of twisted differential operators on $U_e$. Besides, as $\widebar{\mfrak{n}}$ is a nilpotent Lie algebra, the exponential mapping $\exp \colon \widebar{\mfrak{n}} \rarr \widebar{N}$ is an isomorphism of algebraic varieties and induces a canonical isomorphism of algebraic varieties $U_e$ and $\widebar{\mfrak{n}}$. Therefore, the homomorphism $\Phi_\lambda \colon U(\mfrak{g}) \rarr \Gamma(X,\mcal{D}_X^\lambda)$ of $\C$-algebras gives rise to the homomorphism
\begin{align}
  \pi_\lambda \colon U(\mfrak{g}) \rarr \Gamma(X,\mcal{D}_X^\lambda) \rarr \Gamma(U_e,\mcal{D}_X^\lambda) \riso \Gamma(U_e,\mcal{D}_X) \riso \eus{A}_{\widebar{\mfrak{n}}}  \label{eq:pi action abstract}
\end{align}
of $\C$-algebras for $\lambda \in \mfrak{h}^*$, where the Weyl algebra $\eus{A}_{\widebar{\mfrak{n}}}$ of the vector space $\widebar{\mfrak{n}}$ is defined through $\eus{A}_{\widebar{\mfrak{n}}} = \Gamma(\widebar{\mfrak{n}},\mcal{D}_{\widebar{\mfrak{n}}})$.

In addition, there is a nice explicit description of the $\C$-algebra homomorphism \eqref{eq:pi action abstract}. For that reason, let $\{f_\gamma;\, \gamma \in \Delta^+\}$ be a root basis of the opposite nilradical $\widebar{\mfrak{n}}$. We denote by $\{x_\gamma;\, \gamma \in \Delta^+\}$ the linear coordinate functions on $\widebar{\mfrak{n}}$ with respect to the basis $\{f_\gamma;\, \gamma \in \Delta^+\}$ of $\widebar{\mfrak{n}}$. Then the Weyl algebra $\eus{A}_{\widebar{\mfrak{n}}}$ is generated by $\{x_\gamma, \partial_{x_\gamma};\, \gamma \in \Delta^+\}$ together with the canonical commutation relations. The homomorphism $\pi_\lambda \colon U(\mfrak{g}) \rarr \eus{A}_{\widebar{\mfrak{n}}}$ of $\C$-algebras for $\lambda \in \mfrak{h}^*$ is then given by the formula \eqref{eq:pi action general}, see \cite{Krizka-Somberg2017}.
\medskip

Let $\mcal{M}$ be a $i_\alpha^\sharp \mcal{D}_X^\lambda$-module for $\lambda \in \mfrak{h}^*$ and $\alpha \in \Delta^+$. Since $i_\alpha \colon X_\alpha \rarr X$ is a closed embedding, $\mcal{D}_X^\lambda|_{U_e} \simeq \mcal{D}_X|_{U_e}$ as sheaves of rings of twisted differential operators on $U_e$ and $U_e \cap X_\alpha \simeq \mfrak{s}_\alpha^-$, we immediately obtain
\begin{align}
  i_{\alpha*}(\mcal{D}^\lambda_{X \larr X_\alpha} \!\otimes_{i_\alpha^\sharp \mcal{D}_X^\lambda} \mcal{M})|_{U_e} \simeq \C[\partial_{x_\gamma},\gamma \in \Delta^+_\alpha] \otimes_\C i_{\alpha *}(\mcal{M})|_{U_e \simeq \widebar{\mfrak{n}}}, \label{eq:direc image embedding}
\end{align}
where $\Delta^+_\alpha = \Delta^+ \setminus \{\alpha\}$ and the action of the Lie algebra $\mfrak{g}$ on the right hand side is given through the homomorphism $\pi_\lambda \colon U(\mfrak{g}) \rarr \eus{A}_{\widebar{\mfrak{n}}}$ of $\C$-algebras.

If we denote by $M^{\mfrak{s}_\alpha}_{\smash{\mfrak{b} \cap \mfrak{s}_\alpha}}\!((\lambda+\rho-\alpha)_{|\mfrak{h} \cap \mfrak{s}_\alpha})$ the Verma module with highest weight $(\lambda+\rho-\alpha)_{|\mfrak{h} \cap \mfrak{s}_\alpha}$ for the Borel subalgebra $\mfrak{b} \cap \mfrak{s}_\alpha$ and by $N^{\mfrak{s}_\alpha}_{\smash{\mfrak{b} \cap \mfrak{s}_\alpha}}\!((\lambda+\rho)_{|\mfrak{h} \cap \mfrak{s}_\alpha})$ the contragredient Verma module with highest weight $(\lambda+\rho)_{|\mfrak{h} \cap \mfrak{s}_\alpha}$ for the opposite Borel subalgebra $\smash{\widebar{\mfrak{b}}} \cap \mfrak{s}_\alpha$, then we have
\begin{align*}
  M^{\mfrak{s}_\alpha}_{\smash{\mfrak{b} \cap \mfrak{s}_\alpha}}\!((\lambda+\rho-\alpha)_{|\mfrak{h} \cap \mfrak{s}_\alpha}) \simeq \Gamma(X_\alpha,\mcal{L}_{\lambda_\alpha}\!(X_{\alpha,N,e},\mcal{O}_{X_{\alpha,N,e}})) \simeq \C[\partial_{x_\alpha}]
\end{align*}
and
\begin{align*}
   N^{\mfrak{s}_\alpha}_{\smash{\mfrak{b} \cap \mfrak{s}_\alpha}}\!((\lambda+\rho)_{|\mfrak{h} \cap \mfrak{s}_\alpha}) \simeq \Gamma(X_\alpha,\mcal{L}_{\lambda_\alpha}\!(X_{\alpha,\smash{\widebar{N}},e}, \mcal{O}_{X_{\alpha,\smash{\widebar{N}},e}})) \simeq \C[x_\alpha]
\end{align*}
as $\mfrak{s}_\alpha$-modules, where $X_{\alpha,N,e}=N\dot{e}B_\alpha/B_\alpha$ and $X_{\alpha,\smash{\widebar{N}},e} =\widebar{N}\dot{e}B_\alpha/B_\alpha$ are the $N$-orbit and $\widebar{N}$-orbit in $X_\alpha$, respectively. Let us denote by
\begin{align}
  \mcal{K}_\lambda(X_{\alpha,N,e},\mcal{O}_{\alpha,N,e}) = i_{\alpha*}(\mcal{D}^\lambda_{X \larr X_\alpha} \! \otimes_{i_\alpha^\sharp \mcal{D}_X^\lambda}\! \mcal{L}_{\lambda_\alpha}\!(X_{\alpha,N,e},\mcal{O}_{X_{\alpha,N,e}})) \label{eq:K sheaf N-orbit}
\end{align}
and
\begin{align}
  \mcal{K}_\lambda(X_{\alpha,\smash{\widebar{N}},e},\mcal{O}_{\alpha,\smash{\widebar{N}},e}) = i_{\alpha*}(\mcal{D}^\lambda_{X \larr X_\alpha} \! \otimes_{i_\alpha^\sharp \mcal{D}_X^\lambda}\! \mcal{L}_{\lambda_\alpha}\!(X_{\alpha,\smash{\widebar{N}},e}, \mcal{O}_{X_{\alpha,\smash{\widebar{N}},e}})) \label{eq:K sheaf N-orbit op}
\end{align}
for $\lambda \in \mfrak{h}^*$ and $\alpha \in \Delta^+$ coherent $\mcal{D}_X^\lambda$-modules.

The following theorem gives us the geometric realization of Verma modules $M^\mfrak{g}_\mfrak{b}(\lambda)$ and of $\alpha$-Gelfand--Tsetlin modules $W^\mfrak{g}_\mfrak{b}(\lambda,\alpha)$ for $\lambda \in \mfrak{h}^*$ and $\alpha \in \Delta^+$.
\medskip

\theorem{\label{theo-realiz} Let $\lambda \in \mfrak{h}^*$ and $\alpha \in \Delta^+$. Then we have
\begin{align}
  \Gamma(X,\mcal{K}_{\lambda+\rho}(X_{\alpha,N,e}, \mcal{O}_{\alpha,N,e})) \simeq M^\mfrak{g}_\mfrak{b}(\lambda)
\end{align}
and
\begin{align}
  \Gamma(X, \mcal{K}_{\lambda+\rho}(X_{\alpha,\smash{\widebar{N}},e}, \mcal{O}_{\alpha,\smash{\widebar{N}},e})) \simeq W^\mfrak{g}_\mfrak{b}(\lambda,\alpha)
\end{align}
as $\mfrak{g}$-modules.}

\proof{Let $\mcal{M} = \mcal{K}_\lambda(X_{\alpha,N,e}, \mcal{O}_{\alpha,N,e})$. Then we may write
\begin{align*}
  \Gamma(X,\mcal{M}) &\simeq \Gamma(U_e,\mcal{M}) \simeq \C[\partial_{x_\gamma},\, \gamma \in \Delta^+_\alpha] \otimes_\C \Gamma(U_e \cap X_\alpha, \mcal{L}_{\lambda_\alpha}\!(X_{\alpha,N,e},\mcal{O}_{X_{\alpha,N,e}})) \\
  & \simeq \C[\partial_{x_\gamma},\, \gamma \in \Delta^+_\alpha] \otimes_\C \C[\partial_{x_\alpha}] \simeq \C[\partial_{x_\gamma},\, \gamma \in \Delta^+],
\end{align*}
where used \eqref{eq:direc image embedding} and \eqref{eq:K sheaf N-orbit}. Analogously, for $\mcal{M} = \mcal{K}_\lambda(X_{\alpha,\smash{\widebar{N}},e}, \mcal{O}_{\alpha,\smash{\widebar{N}},e})$ we have
\begin{align*}
  \Gamma(X,\mcal{M}) &\simeq \Gamma(U_e,\mcal{M}) \simeq \C[\partial_{x_\gamma},\, \gamma \in \Delta^+_\alpha] \otimes_\C \Gamma(U_e \cap X_\alpha, \mcal{L}_{\lambda_\alpha}\!(X_{\alpha,\smash{\widebar{N}},e}, \mcal{O}_{X_{\alpha,\smash{\widebar{N}},e}})) \\
  & \simeq \C[\partial_{x_\gamma},\, \gamma \in \Delta^+_\alpha] \otimes_\C \C[x_\alpha] \simeq \C[x_\alpha, \partial_{x_\gamma},\, \gamma \in \Delta^+_\alpha].
\end{align*}
By \cite{Krizka-Somberg2017} we have that the $\mfrak{g}$-module $\C[\partial_{x_\alpha},\, \alpha \in \Delta^+]$ is isomorphic to the Verma module $M^\mfrak{g}_\mfrak{b}(\lambda-\rho)$ and by Theorem \ref{thm:Weyl realization} we get that the $\mfrak{g}$-module $\C[x_\alpha,\partial_{x_\gamma},\, \gamma \in \Delta^+_\alpha]$ corresponds to the $\alpha$-Gelfand--Tsetlin module $W^\mfrak{g}_\mfrak{b}(\lambda-\rho,\alpha)$.}

\vspace{-2mm}


\subsection{Geometric realization of $W^\mfrak{g}_\mfrak{p}(\lambda,\alpha)$ -- explicit formulas}

In this section we describe an explicit form of the geometric realization of $\alpha$-Gelfand--Tsetlin modules $W^\mfrak{g}_\mfrak{p}(\lambda,\alpha)$ for $\lambda \in \Lambda^+(\mfrak{p})$ and $\alpha \in \Delta^+_\mfrak{u}$. In the case when $\mfrak{g}$ is a simple Lie algebra, $\mfrak{p}=\mfrak{b}$ and $\alpha$ is the maximal root of $\mfrak{g}$ such realization was constructed in \cite[Theorem 2.14]{Futorny-Krizka2017}.

Let $\{f_\alpha;\, \alpha \in \Delta^+_\mfrak{u}\}$ be a basis of the opposite nilradical $\widebar{\mfrak{u}}$. We denote by $\{x_\alpha;\, \alpha \in \Delta^+_\mfrak{u}\}$ the linear coordinate functions on $\widebar{\mfrak{u}}$ with respect to the basis $\{f_\alpha;\, \alpha \in \Delta^+_\mfrak{u}\}$ of $\widebar{\mfrak{u}}$. Then the Weyl algebra $\eus{A}_{\widebar{\mfrak{u}}}$ of the vector space $\widebar{\mfrak{u}}$ is generated by $\{x_\alpha, \partial_{x_\alpha};\, \alpha \in \Delta^+_\mfrak{u}\}$ together with the canonical commutation relations. For $\lambda \in \Lambda^+(\mfrak{p})$ there is a homomorphism
\begin{align*}
  \pi_\lambda \colon U(\mfrak{g}) \rarr \eus{A}_{\widebar{\mfrak{u}}} \otimes_\C \End \mathbb{F}_{\lambda+\rho_\mfrak{u}}
\end{align*}
of $\C$-algebras given through
\begin{align}
\pi_\lambda(a)= -\sum_{\alpha \in \Delta^+_\mfrak{u}}\bigg[{\ad(u(x))e^{\ad(u(x))} \over e^{\ad(u(x))}-\id}\,(e^{-\ad(u(x))}a)_{\widebar{\mfrak{u}}}\bigg]_\alpha \partial_{x_\alpha} + \sigma_{\lambda+\rho_\mfrak{u}}((e^{-\ad(u(x))}a)_\mfrak{p}) \label{eq:pi action general}
\end{align}
for all $a \in \mfrak{g}$, where the Weyl vector $\rho_\mfrak{u} \in \Lambda^+(\mfrak{p})$ is defined by
\begin{align*}
\rho_\mfrak{u} = {1 \over 2} \sum_{\alpha \in \Delta^+_\mfrak{u}} \alpha,
\end{align*}
$[a]_\alpha$ denotes the $\alpha$-th coordinate of $a \in \widebar{\mfrak{u}}$ with respect to the basis $\{f_\alpha;\, \alpha \in \Delta^+_\mfrak{u}\}$ of $\widebar{\mfrak{u}}$, $a_{\widebar{\mfrak{u}}}$ and $a_\mfrak{p}$ are $\widebar{\mfrak{u}}$-part and $\mfrak{p}$-part of $a \in \mfrak{g}$ with respect to the decomposition $\mfrak{g}=\widebar{\mfrak{u}} \oplus \mfrak{p}$, and finally the element $u(x) \in \C[\widebar{\mfrak{u}}] \otimes_\C \mfrak{g}$ is given by
\begin{align*}
u(x)=\sum_{\alpha \in \Delta^+_\mfrak{u}} x_\alpha f_\alpha.
\end{align*}
Let us note that $\C[\widebar{\mfrak{u}}] \otimes_\C \mfrak{g}$ has the natural structure of a Lie algebra. Hence, we have a well-defined $\C$-linear mapping $\ad(u(x)) \colon \C[\widebar{\mfrak{u}}] \otimes_\C \mfrak{g} \rarr \C[\widebar{\mfrak{u}}] \otimes_\C \mfrak{g}$. In particular, we have
\begin{align}
  \pi_\lambda(a)= - \sum_{\alpha \in\Delta^+_\mfrak{u}} \bigg[{\ad(u(x)) \over e^{\ad(u(x))} - \id}\,a\bigg]_\alpha \partial_{x_\alpha} \label{eq:pi action nilradical}
\end{align}
for $a \in \widebar{\mfrak{u}}$ and
\begin{align}
  \pi_\lambda(a)= \sum_{\alpha \in \Delta^+_\mfrak{u}} [\ad(u(x))a]_\alpha \partial_{x_\alpha} + \sigma_{\lambda+\rho_\mfrak{u}}(a) \label{eq:pi action cartan}
\end{align}
for $a \in \mfrak{l}$.

Let us note that the same construction of the homomorphism $\pi_\lambda \colon U(\mfrak{g}) \rarr \eus{A}_{\widebar{\mfrak{u}}} \otimes_\C \End \mathbb{F}_{\lambda+\rho_\mfrak{u}}$ of $\C$-algebras, as was described in the previous section, works for any parabolic subalgebra $\mfrak{p}$ of $\mfrak{g}$ and any $1$-dimensional $\mfrak{p}$-module $\mathbb{F}_{\lambda+\rho_\mfrak{u}}$, see \cite{Krizka-Somberg2017}. For a general simple finite-dimensional $\mfrak{p}$-module $\mathbb{F}_{\lambda+\rho_\mfrak{u}}$, the expression for $\pi_\lambda$ may be extracted from the formula for the affine Kac--Moody algebra $\widehat{\mfrak{g}}$ associated to $\mfrak{g}$ since $\mfrak{g}$ is a Lie subalgebra of $\widehat{\mfrak{g}}$, see \cite{Futorny-Krizka-Somberg2016}.

The generalized Verma module $M^\mfrak{g}_\mfrak{p}(\lambda)$ for $\lambda \in \Lambda^+(\mfrak{p})$ is realized as $\eus{A}_{\widebar{\mfrak{u}}}/\eus{I}_{\rm V} \otimes_\C \mathbb{F}_{\lambda+2\rho_\mfrak{u}}$, where $\eus{I}_{\rm V}$ is the left ideal of $\eus{A}_{\widebar{\mfrak{u}}}$ defined by $\eus{I}_{\rm V}=(x_\alpha;\, \alpha \in \Delta^+_\mfrak{u})$, see e.g.\ \cite{Krizka-Somberg2017} and \cite{Futorny-Krizka-Somberg2016}. The $\mfrak{g}$-module structure on $\eus{A}_{\widebar{\mfrak{u}}}/\eus{I}_{\rm V} \otimes_\C \mathbb{F}_{\lambda+2\rho_\mfrak{u}}$ is then given through the homomorphism $\pi_{\lambda+\rho_\mfrak{u}} \colon U(\mfrak{g}) \rarr \eus{A}_{\widebar{\mfrak{u}}} \otimes_\C \End \mathbb{F}_{\lambda+2\rho_\mfrak{u}}$ of $\C$-algebras.
\medskip

The $\C$-algebra $\C[\widebar{\mfrak{u}}]$ of polynomial functions on $\widebar{\mfrak{u}}$ has the natural structure of a $Q_+$-graded $\C$-algebra uniquely determined by the requirement $\deg(x_\alpha)=\alpha$ for $\alpha \in \Delta^+_\mfrak{u}$. Therefore, we have a direct sum decomposition
\begin{align*}
  \C[\widebar{\mfrak{u}}] = \bigoplus_{\gamma \in Q_+} \C[\widebar{\mfrak{u}}]_\gamma.
\end{align*}
We say that a differential operator $P \in \eus{A}_{\widebar{\mfrak{u}}}$ has degree $\delta \in Q$ if $P(\C[\widebar{\mfrak{u}}]_\gamma) \subset \C[\widebar{\mfrak{u}}]_{\gamma+\delta}$ for all $\gamma \in Q_+$. The degree of a differential operator defines on the Weyl algebra $\eus{A}_{\widebar{\mfrak{u}}}$ the canonical structure of a $Q$-graded $\C$-algebra, which gives us a direct sum decomposition
\begin{align*}
  \eus{A}_{\widebar{\mfrak{u}}} = \bigoplus_{\gamma \in Q} \eus{A}_{\widebar{\mfrak{u}},\gamma}.
\end{align*}
Further, for $\alpha \in \Delta^+_\mfrak{u}$ we introduce a differential operator $p_\alpha \in \eus{A}_{\widebar{\mfrak{u}}}$ by
\begin{align} \label{eq:p alpha}
  p_\alpha = - \sum_{\gamma \in\Delta^+_\mfrak{u}} \bigg[{\ad(u(x)) \over e^{\ad(u(x))} - \id}\,f_\alpha\bigg]_\gamma \partial_{x_\gamma} = - \sum_{\gamma \in\Delta^+_\mfrak{u}} \partial_{x_\gamma} \bigg[{\ad(u(x)) \over e^{\ad(u(x))} - \id}\,f_\alpha\bigg]_\gamma,
\end{align}
where $u(x) = \sum_{\gamma \in \Delta^+_\mfrak{u}} x_\gamma f_\gamma$. Besides, we may write
\begin{align*}
  p_\alpha = - \partial_{x_\alpha} + q_\alpha
\end{align*}
with
\begin{align} \label{eq:q alpha}
  q_\alpha = - \sum_{\gamma \in\Delta^+_\mfrak{u}} \bigg[\bigg({\ad(u(x)) \over e^{\ad(u(x))} - \id}-\id\bigg)f_\alpha\bigg]_\gamma \partial_{x_\gamma} = - \sum_{\gamma \in\Delta^+_\mfrak{u}} \partial_{x_\gamma} \bigg[\bigg({\ad(u(x)) \over e^{\ad(u(x))} - \id} - \id\bigg)f_\alpha\bigg]_\gamma
\end{align}
for $\alpha\in \Delta^+_\mfrak{u}$. Let us note that the differential operators $p_\alpha$, $q_\alpha$ and $\partial_{x_\alpha}$ have degree $-\alpha$ for $\alpha \in \Delta^+_\mfrak{u}$. We denote by $\eus{N}_\alpha$ the Lie subalgebra of $\eus{A}_{\widebar{\mfrak{u}}}$ generated by the subset $\{\ad(\partial_{x_\alpha})^n(q_\alpha);\, n \in \N_0\}$ and by $\eus{U}_\alpha$ the $\C$-subalgebra of $\eus{A}_{\widebar{\mfrak{u}}}$ generated by $\eus{N}_\alpha$.
\medskip

As $\partial_{x_\alpha}$ and $p_\alpha$ for $\alpha \in \Delta^+_\mfrak{u}$ are locally $\ad$-nilpotent regular elements of $\eus{A}_{\widebar{\mfrak{u}}}$, we may construct left rings of fractions for $\eus{A}_{\widebar{\mfrak{u}}}$ with respect to the multiplicative sets $\{\partial_{x_\alpha}^n;\, n \in \N_0\}$ and $\{p_\alpha^n;\, n\in \N_0\}$.
Let us recall that for $\lambda \in \Lambda^+(\mfrak{p})$ we have an isomorphism
\begin{align}
   M^\mfrak{g}_\mfrak{p}(\lambda-\rho_\mfrak{u}) \simeq \eus{A}_{\widebar{\mfrak{u}}}/\eus{I}_{\rm V} \otimes_\C \mathbb{F}_{\lambda+\rho_\mfrak{u}}
\end{align}
of $U(\mfrak{g})$-modules, which gives rise to an isomorphism $M^\mfrak{g}_\mfrak{p}(\lambda-\rho_\mfrak{u})_{(f_\alpha)} \simeq (\eus{A}_{\widebar{\mfrak{u}}}/\eus{I}_{\rm V} \otimes_\C \mathbb{F}_{\lambda+\rho_\mfrak{u}})_{(f_\alpha)}$ of $U(\mfrak{g})_{(f_\alpha)}$-modules.
Since the $U(\mfrak{g})$-module structure on $\eus{A}_{\widebar{\mfrak{u}}}/\eus{I}_{\rm V} \otimes_\C \mathbb{F}_{\lambda+\rho_\mfrak{u}}$ is given through the homomorphism
\begin{align}
  \pi_\lambda \colon U(\mfrak{g}) \rarr \eus{A}_{\widebar{\mfrak{u}}} \otimes_\C \End \mathbb{F}_{\lambda+\rho_\mfrak{u}} \label{eq:pi homomorphism}
\end{align}
of $\C$-algebras and $\pi_\lambda(f_\alpha) = p_\alpha$ for $\lambda \in \Lambda^+(\mfrak{p})$, we obtain that
\begin{align}
  M^\mfrak{g}_\mfrak{p}(\lambda-\rho_\mfrak{u})_{(f_\alpha)} \simeq (\eus{A}_{\widebar{\mfrak{u}}}/\eus{I}_{\rm V})_{(p_\alpha)} \otimes_\C \mathbb{F}_{\lambda+\rho_\mfrak{u}} \label{eq:localization isomorphism}
\end{align}
as $U(\mfrak{g})_{(f_\alpha)}$-modules, where the $U(\mfrak{g})_{(f_\alpha)}$-module structure on $(\eus{A}_{\widebar{\mfrak{u}}}/\eus{I}_{\rm V})_{(p_\alpha)} \otimes_\C \mathbb{F}_{\lambda+\rho_\mfrak{u}}$ is given by the uniquely determined homomorphism
\begin{align*}
  \pi_{\lambda,\alpha} \colon U(\mfrak{g})_{(f_\alpha)} \rarr (\eus{A}_{\widebar{\mfrak{u}}})_{(p_\alpha)} \otimes_\C \End \mathbb{F}_{\lambda+\rho_\mfrak{u}}
\end{align*}
of $\C$-algebras.

Let us note that the structure of an $\eus{A}_{\widebar{\mfrak{u}}}$-module on $(\eus{A}_{\widebar{\mfrak{u}}}/\eus{I}_{\rm V})_{(p_\alpha)}$ is quite complicated in general. However, we show that
\begin{align*}
  (\eus{A}_{\widebar{\mfrak{u}}}/\eus{I}_{\rm V})_{(p_\alpha)} \simeq (\eus{A}_{\widebar{\mfrak{u}}}/\eus{I}_{\rm V})_{(\partial_{x_\alpha})}
\end{align*}
as modules over the Weyl algebra $\eus{A}_{\widebar{\mfrak{u}}}$.
\medskip

\lemma{\label{lem:locally nilpotent action}
The element $q_\alpha \partial_{x_\alpha}^{-1} \in (\eus{A}_{\widebar{\mfrak{u}}})_{(\partial_{x_\alpha})}$ acts locally nilpotently on $(\eus{A}_{\widebar{\mfrak{u}}}/\eus{I}_{\rm V})_{(\partial_{x_\alpha})}$ for $\alpha \in \Delta^+_\mfrak{u}$.}

\proof{From \eqref{eq:q alpha} we have $\ad(\partial_{x_\alpha})^n(q_\alpha) \in \eus{I}_{\rm V}$ for $n \in \N_0$, which gives us $\eus{N}_\alpha \subset \eus{I}_{\rm V}$ and $\eus{U}_\alpha \subset \eus{I}_{\rm V}$ since $\eus{I}_{\rm V}$ is a left ideal of $\eus{A}_{\widebar{\mfrak{u}}}$. Further, let $\{\eus{F}_k\eus{A}_{\widebar{\mfrak{u}}}\}_{k \in \N_0}$ be a filtration of $\eus{A}_{\widebar{\mfrak{u}}}$ given by the order of a differential operator. Then we have
\begin{align*}
  \eus{N}_\alpha(\eus{F}_k\eus{A}_{\widebar{\mfrak{u}}} + \eus{I}_{\rm V}) \subset \ad(\eus{N}_\alpha)(\eus{F}_k\eus{A}_{\widebar{\mfrak{u}}}) + \eus{I}_{\rm V} \subset \eus{F}_k\eus{A}_{\widebar{\mfrak{u}}} + \eus{I}_{\rm V}
\end{align*}
for $k \in \N_0$. Besides, we have
\begin{align*}
  \eus{N}_\alpha = \bigoplus_{n \in \N_0} \eus{N}_{\alpha,-n\alpha}
\end{align*}
with $\eus{N}_{\alpha,0} =0$ and $\eus{N}_{\alpha,-\alpha} =\C q_\alpha$, which together with the fact that for each $\gamma \in Q$ there exists an integer $n_{k,\gamma} \in \N_0$ such that $\eus{F}_k \eus{A}_{\widebar{\mfrak{u}},\gamma-n\alpha}=0$ for all $n > n_{k,\gamma}$ implies that
\begin{align*}
\eus{U}_{\alpha,-n\alpha}(\eus{F}_k\eus{A}_{\widebar{\mfrak{u}},\gamma} + \eus{I}_{\rm V}) \subset \eus{I}_{\rm V}
\end{align*}
for all $n > n_{k,\gamma}$. In other words, for each $q \in \eus{A}_{\widebar{\mfrak{u}}}/\eus{I}_{\rm V}$ there exists an integer $n_q \in\N_0$ such that $\eus{U}_{\alpha,-n\alpha}q=0$ for all $n>n_q$.

Further, we show by induction on $k$ that
\begin{align*}
  (q_\alpha \partial_{x_\alpha}^{-1})^k = \sum_{j=k}^\infty \partial_{x_\alpha}^{-j} p_{\alpha,k,j},
\end{align*}
where $p_{\alpha,k,j} \in \eus{U}_{\alpha,-j\alpha}$ for $j \geq k$. Moreover, there exists an integer $j_k \in \N_0$ such that $p_{\alpha,k,j}=0$ for $j > j_k$. For $k=0$ we have
\begin{align*}
  p_{\alpha,k,j} =
  \begin{cases}
    1 & \text{if $j=0$},\\
    0 & \text{if $j>0$}.
  \end{cases}
\end{align*}
Let assume that it holds for some $k \in \N_0$. Then we may write
\begin{align*}
  (q_\alpha \partial_{x_\alpha}^{-1})^{k+1} &= q_\alpha \partial_{x_\alpha}^{-1} \sum_{k=j}^\infty \partial_{x_\alpha}^{-j} p_{\alpha,k,j} = \sum_{k=j}^\infty \partial_{x_\alpha}^{-j-1} \partial_{x_\alpha}^{j+1} q_\alpha \partial_{x_\alpha}^{-j-1} p_{\alpha,k,j} \\
  &= \sum_{j=k}^\infty \partial_{x_\alpha}^{-j-1} \sum_{\ell=0}^\infty \binom{j+\ell}{\ell} \partial_{x_\alpha}^{-\ell} \ad(\partial_{x_\alpha})^\ell(q_\alpha) p_{\alpha,k,j} \\
  &= \sum_{j=k+1}^\infty \partial_{x_\alpha}^{-j} \sum_{\ell=0}^{j-k-1} \binom{j-1}{\ell} \ad(\partial_{x_\alpha})^\ell(q_\alpha)p_{\alpha,k,j-\ell-1},
\end{align*}
where we used
\begin{align*}
  \partial_{x_\alpha}^n q \partial_{x_\alpha}^{-n} = \sum_{k=0}^\infty \binom{n+k-1}{k} \partial_{x_\alpha}^{-k} \ad(\partial_{x_\alpha})^k(q)
\end{align*}
for $n \in \N_0$ and $q \in \eus{A}_{\widebar{\mfrak{u}}}$ in the third equality,
which gives us
\begin{align*}
  p_{\alpha,k+1,j}= \sum_{\ell=0}^{j-k-1} \binom{j-1}{\ell} \ad(\partial_{x_\alpha})^\ell(q_\alpha)p_{\alpha,k,j-\ell-1}
\end{align*}
for $j \geq k+1$. Therefore, we have $p_{\alpha,k+1,j} \in \eus{U}_{\alpha,-j\alpha}$ for $j\geq k+1$ and $j_{k+1} = j_k + m +1$, where $m \in \N_0$ satisfies $\ad(\partial_{x_\alpha})^m(q_\alpha)=0$.

Furthermore, for $n\in \N_0$ and $q \in \eus{A}_{\widebar{\mfrak{u}}}$ we have
\begin{align*}
  p_{\alpha,k,j}\partial_{x_\alpha}^{-n}q &= \partial_{x_\alpha}^{-n} \partial_{x_\alpha}^n p_{\alpha,k,j} \partial_{x_\alpha}^{-n} q = \partial_{x_\alpha}^{-n} \sum_{\ell=0}^\infty \binom{n+\ell-1}{\ell} \partial_{x_\alpha}^{-\ell} \ad(\partial_{x_\alpha})^\ell(p_{\alpha,k,j})q
\end{align*}
for $j,k \in \N_0$ and $j\geq k$. As we have $p_{\alpha,k,j}\in \eus{U}_{\alpha,-j\alpha}$, we obtain $\ad(\partial_{x_\alpha})^\ell(p_{\alpha,k,j}) \in \eus{U}_{\alpha,-(\ell+j)\alpha}$ for $\ell \in \N_0$. Since for $q \in \eus{A}_{\widebar{\mfrak{u}}}/\eus{I}_{\rm V}$ there exists an integer $k_q \in\N_0$ such that $\eus{U}_{\alpha,-k\alpha}q=0$ for all $k>k_q$, we get that $p_{\alpha,k,j}\partial_{x_\alpha}^{-n}q=0$ for all $k>k_q$ and $j \geq k$. Therefore, we have $(q_\alpha \partial_{x_\alpha}^{-1})^k \partial_{x_\alpha}^{-n}q = 0$ for $k > k_q$. Hence, we have the required statement.}

\lemma{\label{lem:p_alpha inverse}
For $\alpha \in \Delta^+_\mfrak{u}$ the $\C$-linear mapping
\begin{align*}
  \varphi_\alpha \colon (\eus{A}_{\widebar{\mfrak{u}}}/\eus{I}_{\rm V})_{(\partial_{x_\alpha})} \rarr (\eus{A}_{\widebar{\mfrak{u}}}/\eus{I}_{\rm V})_{(\partial_{x_\alpha})}
\end{align*}
defined by
\begin{align}
  \varphi_\alpha(a)= - \partial_{x_\alpha}^{-1} \sum_{k=0}^\infty (q_\alpha \partial_{x_\alpha}^{-1})^k a \label{eq:p_alpha inverse}
\end{align}
for $a \in (\eus{A}_{\widebar{\mfrak{u}}}/\eus{I}_{\rm V})_{(\partial_{x_\alpha})}$ satisfies the relations
\begin{align*}
  p_\alpha \varphi_\alpha(a) = a \qquad \text{and} \qquad \varphi_\alpha(p_\alpha a) = a
\end{align*}
for $a \in (\eus{A}_{\widebar{\mfrak{u}}}/\eus{I}_{\rm V})_{(\partial_{x_\alpha})}$.}

\proof{By Lemma \ref{lem:locally nilpotent action} we have that the $\C$-linear mapping $\varphi_\alpha$ is well defined since the sum on the right hand side is finite. Further, we may write
\begin{align*}
  p_\alpha \varphi_\alpha(a) = (\partial_{x_\alpha} - q_\alpha) \partial_{x_\alpha}^{-1} \sum_{k=0}^\infty (q_\alpha \partial_{x_\alpha}^{-1})^k a = \sum_{k=0}^\infty (q_\alpha \partial_{x_\alpha}^{-1})^k a - \sum_{k=0}^\infty (q_\alpha \partial_{x_\alpha}^{-1})^{k+1} a = a
\end{align*}
and
\begin{align*}
  \varphi_\alpha(p_\alpha a) &= \partial_{x_\alpha}^{-1} \sum_{k=0}^\infty (q_\alpha \partial_{x_\alpha}^{-1})^k (\partial_{x_\alpha} - q_\alpha) a = \partial_{x_\alpha}^{-1} \sum_{k=0}^\infty (q_\alpha \partial_{x_\alpha}^{-1})^k (1 - q_\alpha\partial_{x_\alpha}^{-1})\partial_{x_\alpha} a \\
  &= \partial_{x_\alpha}^{-1} \bigg(\sum_{k=0}^\infty (q_\alpha \partial_{x_\alpha}^{-1})^k  - \sum_{k=0}^\infty (q_\alpha \partial_{x_\alpha}^{-1})^{k+1}\!\bigg)\partial_{x_\alpha} a =a
\end{align*}
for all $a \in (\eus{A}_{\widebar{\mfrak{u}}}/\eus{I}_{\rm V})_{(\partial_{x_\alpha})}$.}

\lemma{Let $\alpha \in \Delta^+_\mfrak{u}$. Then for all $a \in (\eus{A}_{\widebar{\mfrak{u}}}/\eus{I}_{\rm V})_{(\partial_{x_\alpha})}$ we have
\begin{enumerate}[topsep=3pt,itemsep=0pt]
  \item[i)] $\varphi_\alpha^n(p_\alpha^ma) = p_\alpha^m \varphi_\alpha^n(a) = \varphi_\alpha^{n-m}(a)$ for $n,m \in \N_0$ and $n \geq m$;
  \item[ii)] $p_\alpha^nq\varphi_\alpha^n(a) = \sum_{k=0}^\infty \binom{n+k-1}{k} \varphi_\alpha^k(\ad(p_\alpha)^k(q)a)$ for $n \in \N_0$ and $q \in \eus{A}_{\widebar{\mfrak{u}}}$.
\end{enumerate}}

\proof{The first statement follows immediately from Lemma \ref{lem:p_alpha inverse}. Further, for $n\in \N_0$ and $q \in \eus{A}_{\widebar{\mfrak{u}}}$ we have
\begin{align*}
  p_\alpha q\varphi_\alpha(a) = \sum_{k=0}^n \varphi_\alpha^k(\ad(p_\alpha)^k(q)a) + \varphi_\alpha^n(\ad(p_\alpha)^{n+1}(q)\varphi_\alpha(a))
\end{align*}
for $a \in (\eus{A}_{\widebar{\mfrak{u}}}/\eus{I}_{\rm V})_{(\partial_{x_\alpha})}$, which follows by induction on $n$. Indeed, for $n=0$ we have
\begin{align*}
  p_\alpha q\varphi_\alpha(a) = qp_\alpha\varphi_\alpha(a) + \ad(p_\alpha)(q)\varphi_\alpha(a) = qa + \ad(p_\alpha)(q)\varphi_\alpha(a)
\end{align*}
for $a \in (\eus{A}_{\widebar{\mfrak{u}}}/\eus{I}_{\rm V})_{(\partial_{x_\alpha})}$. Let us assume that it holds for some $n \in \N_0$. Then we may write
\begin{align*}
  p_\alpha q\varphi_\alpha(a) &= \sum_{k=0}^n \varphi_\alpha^k(\ad(p_\alpha)^k(q)a) + \varphi_\alpha^n(\ad(p_\alpha)^{n+1}(q)\varphi_\alpha(a)) \\
  &= \sum_{k=0}^n \varphi_\alpha^k(\ad(p_\alpha)^k(q)a) + \varphi_\alpha^{n+1}(p_\alpha\ad(p_\alpha)^{n+1}(q)\varphi_\alpha(a)) \\
  &= \sum_{k=0}^n \varphi_\alpha^k(\ad(p_\alpha)^k(q)a) + \varphi_\alpha^{n+1}(\ad(p_\alpha)^{n+1}(q)a) + \varphi_\alpha^{n+1}(\ad(p_\alpha)^{n+2}(q)\varphi_\alpha(a))
\end{align*}
for $a \in (\eus{A}_{\widebar{\mfrak{u}}}/\eus{I}_{\rm V})_{(\partial_{x_\alpha})}$. As there exists an integer $n_q \in \N_0$ such that $\ad(p_\alpha)^n(q) =0$ for all $n>n_q$, we obtain
\begin{align*}
  p_\alpha q \varphi_\alpha(a) = \sum_{k=0}^\infty \varphi_\alpha^k(\ad(p_\alpha)^k(q)a)
\end{align*}
for $a \in (\eus{A}_{\widebar{\mfrak{u}}}/\eus{I}_{\rm V})_{(\partial_{x_\alpha})}$. Therefore, we proved the second statement for $n=1$. The rest of the proof is by induction on $n$.
Let us assume that it holds for some $n \in \N_0$. Then we may write
\begin{align*}
  p_\alpha^{n+1} q \varphi_\alpha^{n+1}(a) &= p_\alpha \sum_{k=0}^\infty \binom{n+k-1}{k} \varphi_\alpha^k(\ad(p_\alpha)^k(q)\varphi_\alpha(a)) \\
  &= \sum_{k=0}^\infty \binom{n+k-1}{k} \varphi_\alpha^k(p_\alpha\ad(p_\alpha)^k(q)\varphi_\alpha(a)) \\
  &= \sum_{k=0}^\infty \binom{n+k-1}{k} \sum_{\ell=0}^\infty \varphi_\alpha^{k+\ell}( \ad(p_\alpha)^{k+\ell}(q)a) \\
  &= \sum_{k=0}^\infty \sum_{j=0}^k \binom{n+j-1}{j} \varphi_\alpha^k(\ad(p_\alpha)^k(q)a) = \sum_{k=0}^\infty \binom{n+k}{k} \varphi_\alpha^k(\ad(p_\alpha)^k(q)a)
\end{align*}
for $a \in (\eus{A}_{\widebar{\mfrak{u}}}/\eus{I}_{\rm V})_{(\partial_{x_\alpha})}$.}

The previous lemma enables us to define the structure of an $(\eus{A}_{\widebar{\mfrak{u}}})_{(p_\alpha)}$-module on $(\eus{A}_{\widebar{\mfrak{u}}}/\eus{I}_{\rm V})_{(\partial_{x_\alpha})}$ for $\alpha \in \Delta^+_\mfrak{u}$ by
\begin{align*}
  p_\alpha^{-n}qa = \varphi_\alpha^n(qa),
\end{align*}
where $n \in \N_0$, $q \in \eus{A}_{\widebar{\mfrak{u}}}$ and $a \in (\eus{A}_{\widebar{\mfrak{u}}}/\eus{I}_{\rm V})_{(\partial_{x_\alpha})}$.
\medskip

\proposition{\label{prop:Weyl algebra modules isomorphism}
Let $\alpha \in \Delta^+_\mfrak{u}$. Then the $\C$-linear mapping
\begin{align*}
  \Phi_\alpha \colon (\eus{A}_{\widebar{\mfrak{u}}}/\eus{I}_{\rm V})_{(p_\alpha)} \rarr (\eus{A}_{\widebar{\mfrak{u}}}/\eus{I}_{\rm V})_{(\partial_{x_\alpha})}
\end{align*}
defined through
\begin{align}
  \Phi_\alpha(p_\alpha^{-n}a) = \varphi_\alpha^n(a)
\end{align}
for $n \in \N_0$ and $a \in \eus{A}_{\widebar{\mfrak{u}}}/\eus{I}_{\rm V}$ is an isomorphism
of $(\eus{A}_{\widebar{\mfrak{u}}})_{(p_\alpha)}$-modules.}

\proof{Since $\Phi_\alpha$ is a homomorphism of $(\eus{A}_{\widebar{\mfrak{u}}})_{(p_\alpha)}$-modules based on the previous considerations,
we only need to prove that the $\C$-linear mapping $\Phi_\alpha$ is injective and surjective. We may write
\begin{align*}
\Phi_\alpha({\textstyle \prod_{\gamma \in \Delta^+_\mfrak{u}}} \partial_{x_\gamma}^{n_\gamma} p_\alpha^{-n}1) = {\textstyle \prod_{\gamma \in \Delta^+_\mfrak{u}}} \partial_{x_\gamma}^{n_\gamma} \varphi_\alpha^n(1) = (-1)^n {\textstyle \prod_{\gamma \in \Delta^+_\mfrak{u}}} \partial_{x_\gamma}^{n_\gamma}  \partial_{x_\alpha}^{-n}
\end{align*}
for $n \in \N_0$ and $n_\gamma \in \N_0$ for $\gamma \in \Delta^+_\mfrak{u}$, which implies that $\Phi_\alpha$ is surjective. To prove the injectivity of $\Phi_\alpha$ let us assume that $\Phi_\alpha(a) = 0$ for some $a \in (\eus{A}_{\widebar{\mfrak{u}}}/\eus{I}_{\rm V})_{(p_\alpha)}$. Then there exists an integer $n \in \N_0$ such that $p_\alpha^n a \in \eus{A}_{\widebar{\mfrak{u}}}/\eus{I}_{\rm V} \subset (\eus{A}_{\widebar{\mfrak{u}}}/\eus{I}_{\rm V})_{(p_\alpha)}$. Hence, we have
\begin{align*}
 0 = p_\alpha^n \Phi_\alpha(a) = \Phi_\alpha(p_\alpha^na) = p_\alpha^n a \in \eus{A}_{\widebar{\mfrak{u}}}/\eus{I}_{\rm V} \subset (\eus{A}_{\widebar{\mfrak{u}}}/\eus{I}_{\rm V})_{(\partial_{x_\alpha})},
\end{align*}
which gives us $a=0$. Therefore, we proved that $\Phi_\alpha$ is an isomorphism of $(\eus{A}_{\widebar{\mfrak{u}}})_{(p_\alpha)}$-modules.}

Let us recall that for $\lambda \in \Lambda^+(\mfrak{p})$ and $\alpha \in \Delta^+_\mfrak{u}$ we have an isomorphism
\begin{align*}
  M^\mfrak{g}_\mfrak{p}(\lambda-\rho_\mfrak{u})_{(f_\alpha)} \simeq (\eus{A}_{\widebar{\mfrak{u}}}/\eus{I}_{\rm V})_{(p_\alpha)} \otimes_\C \mathbb{F}_{\lambda+\rho_\mfrak{u}}
\end{align*}
of $U(\mfrak{g})$-modules, where the $U(\mfrak{g})$-module structure on $(\eus{A}_{\widebar{\mfrak{u}}}/\eus{I}_{\rm V})_{(p_\alpha)} \otimes_\C \mathbb{F}_{\lambda+\rho_\mfrak{u}}$ is given through the homomorphism \eqref{eq:pi homomorphism} of $\C$-algebras. By Proposition \ref{prop:Weyl algebra modules isomorphism} we have an isomorphism $(\eus{A}_{\widebar{\mfrak{u}}}/\eus{I}_{\rm V})_{(p_\alpha)} \simeq (\eus{A}_{\widebar{\mfrak{u}}}/\eus{I}_{\rm V})_{(\partial_{x_\alpha})}$ of $\eus{A}_{\widebar{\mfrak{u}}}$-modules, which gives rise to an isomorphism
\begin{align}
  M^\mfrak{g}_\mfrak{p}(\lambda-\rho_\mfrak{u})_{(f_\alpha)} \simeq (\eus{A}_{\widebar{\mfrak{u}}}/\eus{I}_{\rm V})_{(\partial_{x_\alpha})} \otimes_\C \mathbb{F}_{\lambda+\rho_\mfrak{u}}  \label{eq:Verma module localization}
\end{align}
of $U(\mfrak{g})$-modules, where the action of $U(\mfrak{g})$ on $(\eus{A}_{\widebar{\mfrak{u}}}/\eus{I}_{\rm V})_{(\partial_{x_\alpha})} \otimes_\C \mathbb{F}_{\lambda+\rho_\mfrak{u}}$ is given by the homomorphism \eqref{eq:pi homomorphism} of $\C$-algebras.
\medskip

We denote by $\Delta^+_{\mfrak{u},\alpha}$ the set $\Delta^+_\mfrak{u} \setminus \{\alpha\}$ for $\alpha \in \Delta^+_\mfrak{u}$. Further, let us consider $\eus{A}_{\widebar{\mfrak{u}}}$-modules
\begin{align}
  \eus{F}_{\widebar{\mfrak{u}}} = \C[\partial_{x_\gamma},\gamma \in \Delta^+_\mfrak{u}] \qquad \text{and} \qquad
  \eus{F}_{\widebar{\mfrak{u}},\alpha} = \C[x_\alpha,\partial_{x_\gamma},\gamma \in \Delta^+_{\mfrak{u},\alpha}]
\end{align}
for $\alpha \in \Delta^+_\mfrak{u}$, where $\eus{F}_{\widebar{\mfrak{u}}}$ and $\eus{F}_{\widebar{\mfrak{u}},\alpha}$ are endowed with the structure of $\eus{A}_{\widebar{\mfrak{u}}}$-modules by means of the canonical isomorphisms of vector spaces
\begin{align}
  \eus{A}_{\widebar{\mfrak{u}}}/\eus{I}_{\rm V} \simeq \eus{F}_{\widebar{\mfrak{u}}} \qquad \text{and} \qquad
  \eus{A}_{\widebar{\mfrak{u}}}/\eus{I}_{{\rm GT},\alpha} \simeq \eus{F}_{\widebar{\mfrak{u}},\alpha}.  \label{eq:A-module structure}
\end{align}
The left ideals $\eus{I}_{\rm V}$ and $\eus{I}_{{\rm GT},\alpha}$ of the Weyl algebra $\eus{A}_{\widebar{\mfrak{u}}}$ are defined through $\eus{I}_{\rm V}=(x_\gamma,\gamma \in \Delta^+_\mfrak{u})$ and $\eus{I}_{{\rm GT},\alpha}=(\partial_{x_\alpha}, x_\gamma,\gamma \in \Delta^+_{\mfrak{u},\alpha})$. In addition, we have a short exact sequence
\begin{align}
  0 \rarr \eus{F}_{\widebar{\mfrak{u}}} \rarr (\eus{F}_{\widebar{\mfrak{u}}})_{(\partial_{x_\alpha})} \rarr \eus{F}_{\widebar{\mfrak{u}},\alpha} \rarr 0 \label{eq:ses fock}
\end{align}
of $\eus{A}_{\widebar{\mfrak{u}}}$-modules, where the surjective homomorphism of $\eus{A}_{\widebar{\mfrak{u}}}$-modules from $(\eus{F}_{\widebar{\mfrak{u}}})_{(\partial_{x_\alpha})}$ to $\eus{F}_{\widebar{\mfrak{u}},\alpha}$ is given by
\begin{align*}
  {\textstyle \prod_{\gamma \in \Delta^+_\mfrak{u}}}\, \partial_{x_\gamma}^{n_\gamma} \mapsto
  \begin{cases}
    {x_\alpha^{-n_\alpha-1} \over (-n_\alpha-1)!}  \prod_{\gamma \in \Delta^+_{\mfrak{u},\alpha}} \partial_{x_\gamma}^{n_\gamma}  & \text{if $n_\alpha < 0$}, \\
    0 & \text{if $n_\alpha \geq 0$}
  \end{cases}
\end{align*}
for $n_\alpha \in \Z$ and $n_\gamma \in \N_0$, $\gamma \in \Delta^+_{\mfrak{u},\alpha}$.
\medskip

\theorem{\label{thm:Weyl realization}
Let $\lambda \in \Lambda^+(\mfrak{p})$ and $\alpha \in \Delta^+_\mfrak{u}$. Then we have
\begin{align}
  M^\mfrak{g}_\mfrak{p}(\lambda) \simeq \eus{A}_{\widebar{\mfrak{u}}}/\eus{I}_{\rm V} \otimes_\C \mathbb{F}_{\lambda+2\rho_\mfrak{u}} \simeq \eus{F}_{\widebar{\mfrak{u}}} \otimes_\C \mathbb{F}_{\lambda+2\rho_\mfrak{u}}
\end{align}
and
\begin{align}
  W^\mfrak{g}_\mfrak{p}(\lambda,\alpha) \simeq \eus{A}_{\widebar{\mfrak{u}}}/\eus{I}_{{\rm GT},\alpha} \otimes_\C \mathbb{F}_{\lambda+2\rho_\mfrak{u}} \simeq \eus{F}_{\widebar{\mfrak{u}},\alpha} \otimes_\C \mathbb{F}_{\lambda+2\rho_\mfrak{u}},
\end{align}
where $\eus{I}_{\rm V}$ and $\eus{I}_{{\rm GT},\alpha}$ are the left ideals of the Weyl algebra $\eus{A}_{\widebar{\mfrak{u}}}$ defined by $\eus{I}_{\rm V}=(x_\gamma,\gamma \in \Delta^+_\mfrak{u})$ and $\eus{I}_{{\rm GT},\alpha}=(\partial_{x_\alpha}, x_\gamma,\gamma \in \Delta^+_{\mfrak{u},\alpha})$.}

\proof{From the previous considerations and definition of $W^\mfrak{g}_\mfrak{p}(\lambda,\alpha)$ for $\lambda \in \Lambda^+(\mfrak{p})$ and $\alpha \in \Delta^+_\mfrak{u}$ we have
\begin{align*}
  W^\mfrak{g}_\mfrak{p}(\lambda,\alpha) \simeq M^\mfrak{g}_\mfrak{p}(\lambda)_{(f_\alpha)}/M^\mfrak{g}_\mfrak{p}(\lambda) \simeq
  ((\eus{F}_{\widebar{\mfrak{u}}})_{(\partial_{x_\alpha})}/\eus{F}_{\widebar{\mfrak{u}}}) \otimes_\C \mathbb{F}_{\lambda+2\rho_\mfrak{u}} \simeq \eus{F}_{\widebar{\mfrak{u}},\alpha} \otimes_\C \mathbb{F}_{\lambda+2\rho_\mfrak{u}},
\end{align*}
where we used \eqref{eq:Verma module localization}, \eqref{eq:A-module structure} and \eqref{eq:ses fock}.}

\vspace{-2mm}


\section*{Acknowledgments}

V.\,F.\ is supported in part by the CNPq (304467/2017-0 and 200783/2018-1); L.\,K.\ gratefully acknowledges the hospitality and excellent working conditions at the University of Sa\~{o} Paulo where this work was done. V.\,F.\ gratefully acknowledges the hospitality and excellent working conditions at the University of California Berkeley where this work was completed.



\begin{thebibliography}{10}

\bibitem{Andersen-Lauritzen2003}
Henning~H. Andersen and Niels Lauritzen, \emph{{Twisted Verma modules}},
  Studies in Memory of Issai Schur, Progress in Mathematics, vol. 210,
  Birkhäuser, Boston, 2003, pp.~1--26.

\bibitem{Andersen-Stroppel2003}
Henning~H. Andersen and Catharina Stroppel, \emph{{Twisting functors on
  $\mathcal{O}$}}, Represent. Theory \textbf{7} (2003), 681--699.

\bibitem{Arkhipov1997}
Sergey~M. Arkhipov, \emph{{Semi-infinite cohomology of associative algebras and
  bar duality}}, Internat. Math. Res. Notices (1997), no.~17, 833--863.

\bibitem{Arkhipov2004}
\bysame, \emph{{Algebraic construction of contragradient quasi-Verma modules in
  positive characteristic}}, Adv. Stud. Pure Math. \textbf{40} (2004), 27--68.

\bibitem{Beilinson-Bernstein1981}
Alexander~A. Beilinson and Joseph~N. Bernstein, \emph{{Localisation de
  $\mathfrak{g}$-modules}}, C. R. Acad. Sci. Paris Sér. I Math. \textbf{292}
  (1981), 15--18.


\bibitem{Deodhar1980}
Vinay~V. Deodhar, \emph{{On a construction of representations and a problem of
  Enright}}, Invent. Math. \textbf{57} (1980), no.~2, 101--118.


\bibitem{Drozd-Futorny-Ovsienko1994}
Yuri~A. Drozd, Vyacheslav Futorny, and Serge Ovsienko, \emph{{Harish-Chandra
  subalgebras and Gelfand-Zetlin modules}}, Finite-dimensional algebras and
  related topics (Ottawa, ON, 1992), NATO Adv. Sci. Inst. Ser. C Math. Phys.
  Sci., vol. 424, Kluwer Acad. Publ., Dordrecht, 1994, pp.~79--93.

\bibitem{Drozd-Ovsienko-Futorny1991}
Yuri~A. Drozd, Serge Ovsienko, and Vyacheslav Futorny, \emph{{On Gelfand-Zetlin
  modules}}, Rend. Circ. Mat. Palermo (2) Suppl. \textbf{26} (1991), 143--147.

\bibitem{Early-Mazorchuk-Vishnyakova2017}
Nick Early, Volodymyr Mazorchuk, and Elizaveta Vishnyakova, \emph{{Canonical
  Gelfand-Zeitlin modules over orthogonal Gelfand-Zeitlin algebras}}, {\tt
  arXiv:1709.01553} (2017).

\bibitem{Futorny-Grantcharov-Ramirez2014}
Vyacheslav Futorny, Dimitar Grantcharov, and Luis~E. Ramírez, \emph{{On the
  classification of irreducible Gelfand--Tsetlin modules of
  $\mathfrak{sl}(3)$}}, Recent Advances in Representation Theory, Quantum
  Groups, Algebraic Geometry, and Related Topics, Contemporary Mathematics,
  vol. 623, Amer. Math. Soc., Providence, RI, 2014, pp.~63--79.

\bibitem{Futorny-Grantcharov-Ramirez2016}
\bysame, \emph{{Singular Gelfand-Tsetlin modules of $\mathfrak{gl}(n)$}}, Adv.
  Math. \textbf{290} (2016), 453--482.

\bibitem{Futorny-Grantcharov-Ramirez2017a}
\bysame, \emph{{Drinfeld category and the classification of singular
  Gelfand-Tsetlin $\mathfrak{gl}_n$-modules}}, {\tt arXiv:1704.01209} (2017).

\bibitem{Futorny-Grantcharov-Ramirez2017}
\bysame, \emph{{New singular Gelfand-Tsetlin $\mathfrak{gl}(n)$-modules of
  index $2$}}, Comm. Math. Phys. \textbf{355} (2017), no.~3, 1209--1241.

\bibitem{FGR-sl(3)}
\bysame, \emph{{Classification of irreducible Gelfand-Tsetlin modules of $sl(3)$}},
2018, arXiv:1812.07137.

\bibitem{Futorny-Grantcharov-Ramirez-Zadunaisky2018a}
Vyacheslav Futorny, Dimitar Grantcharov, Luis~E. Ramírez, and Pablo Zadunaisky,
  \emph{{Bounds of Gelfand-Tsetlin multiplicities and tableaux realizations of
  Verma modules}}, {\tt arXiv:1811.07992} (2018).

\bibitem{Futorny-Grantcharov-Ramirez-Zadunaisky2018}
\bysame, \emph{{Gelfand-Tsetlin theory for rational Galois algebras}}, {\tt
  arXiv:1801.09316} (2018).

\bibitem{Futorny-Krizka2017}
Vyacheslav Futorny and Libor Křižka, \emph{{Geometric construction of
  Gelfand--Tsetlin modules over simple Lie algebras}}, {\tt arXiv:1712.03700}
  (2017) (to appear in J. Pure Appl. Algebra).

\bibitem{Futorny-Krizka-Somberg2016}
Vyacheslav Futorny, Libor Křižka, and Petr Somberg, \emph{{Geometric
  realizations of affine Kac--Moody algebras}}, {\tt arXiv:1610.07973} (2016)
  (to appear in J. Algebra).

\bibitem{Futorny-Ovsienko-Saorin2011}
Vyacheslav Futorny, Serge Ovsienko, and Manuel Saorín, \emph{{Torsion theories
  induced from commutative subalgebras}}, J. Pure Appl. Algebra \textbf{215}
  (2011), no.~12, 2937--2948.

\bibitem{Futorny-Ramirez-Zhang2016}
Vyacheslav Futorny, Luis~E. Ramírez, and Jian Zhang, \emph{{Combinatorial
  construction of Gelfand-Tsetlin modules for $\mathfrak{gl}_n$}}, {\tt
  arXiv:1611.07908} (2016) (to appear in Adv. Math).

\bibitem{Graev2007}
Mark~I. Graev, \emph{{A continuous analogue of Gelfand-Tsetlin schemes and a
  realization of the principal series of irreducible unitary representations of
  the group ${\rm GL}(n,\mathbb{C})$ in the space of functions on the manifold
  of these schemes}}, Dokl. Akad. Nauk. \textbf{412} (2007), no.~2, 154--158.

\bibitem{Guillemin-Sternberg1983}
Victor Guillemin and Shlomo Sternberg, \emph{{The Gelfand-Cetlin system and
  quantization of the complex flag manifolds}}, J. Funct. Anal. \textbf{52}
  (1983), no.~1, 106--128.

\bibitem{Hartwig2017}
Jonas~T. Hartwig, \emph{{Principal Galois orders and Gelfand-Zeitlin modules}},
  {\tt arXiv:1710.04186} (2017).

\bibitem{Kamnitzer-Tingley-Webster-Weekes-Yacobi2018}
Joel Kamnitzer, Peter Tingley, Ben Webster, Alex Weekes, and Oded Yacobi,
  \emph{{On category $\mathcal{O}$ for affine Grassmannian slices and
  categorified tensor products}}, {\tt arXiv:1806.07519} (2018).

\bibitem{Kashiwara1989}
Masaki Kashiwara, \emph{{Representaion theory and $\mathcal{D}$-modules on flag
  varieties}}, Astérisque \textbf{173-174} (1989), 55--109.

\bibitem{Krizka-Somberg2015b}
Libor Křižka and Petr Somberg, \emph{{On the composition structure of the
  twisted Verma modules for $\mathfrak{sl}(3,\mathbb{C})$}}, Arch. Math.
  \textbf{51} (2015), no.~5, 315--329.

\bibitem{Krizka-Somberg2017}
\bysame, \emph{{Algebraic analysis of scalar generalized Verma modules of
  Heisenberg parabolic type I: $A_n$-series}}, Transform. Groups \textbf{22}
  (2017), no.~2, 403--451.

\bibitem{Kleshchev2018}
Alexander Kleshchev, 2018, private communication.

\bibitem{Kostant-Wallach2006}
Bertram Kostant and Nolan Wallach, \emph{{Gelfand-Zeitlin theory from the
  perspective of classical mechanics. I}}, Studies in Lie theory, Progress in
  Mathematics, vol. 243, Birkhäuser Boston, Boston, MA, 2006, pp.~319--364.

\bibitem{Mathieu2000}
Olivier Mathieu, \emph{{Classification of irreducible weight modules}}, Ann.
  Inst. Fourier (Grenoble) \textbf{50} (2000), no.~2, 537--592.

\bibitem{Mazorchuk-Vishnyakova2018}
Volodymyr Mazorchuk and Elizaveta Vishnyakova, \emph{{Harish-Chandra modules
  over invariant subalgebras in a skew-group ring}}, {\tt arXiv:1811.00332}
  (2018).

\bibitem{Musson2019}
Ian~M. Musson, \emph{{Twisting functors and generalized Verma modules}}, Proc.
  Amer. Math. Soc. \textbf{147} (2019), no.~3, 1013--1022.

\bibitem{Ramirez-Zadunaisky2018}
Luis~E. Ramírez and Pablo Zadunaisky, \emph{{Gelfand-Tsetlin modules over
  $\mathfrak{gl}(n,\mathbb{C})$ with arbitrary characters}}, J. Algebra
  \textbf{502} (2018), 328--346.

\bibitem{Vishnyakova2017}
Elizaveta Vishnyakova, \emph{{Geometric approach to $p$-singular
  Gelfand-Tsetlin $\mathfrak{gl}_n$-modules}}, {\tt arXiv:1705.05793} (2017).

\bibitem{Vishnyakova2018}
\bysame, \emph{{A geometric approach to $1$-singular Gelfand-Tsetlin
  $\mathfrak{gl}_n$-modules}}, Differential Geom. Appl. \textbf{56} (2018),
  155--160.

\bibitem{Zadunaisky2017}
Pablo Zadunaisky, \emph{{A new way to construct $1$-singular Gelfand-Tsetlin
  modules}}, Algebra Discrete Math. \textbf{23} (2017), no.~1, 180--193.

\end{thebibliography}

\providecommand{\bysame}{\leavevmode\hbox to3em{\hrulefill}\thinspace}
\providecommand{\MR}{\relax\ifhmode\unskip\space\fi MR }
\providecommand{\MRhref}[2]{%
  \href{http://www.ams.org/mathscinet-getitem?mr=#1}{#2}
}
\providecommand{\href}[2]{#2}

\end{document}